\documentclass{article}
\usepackage[colorlinks=true,citecolor=blue,urlcolor=blue,linkcolor=blue]{hyperref}
\usepackage[margin=1in]{geometry}
\usepackage[utf8]{inputenc}
\usepackage[T1]{fontenc}
\usepackage{url}
\usepackage{booktabs}
\usepackage{amsfonts}
\usepackage{nicefrac}
\usepackage{microtype}
\usepackage{xcolor}
\usepackage{grffile}
\usepackage[shortlabels]{enumitem}
\usepackage[round]{natbib}

\usepackage{tocvsec2}

\usepackage{amsmath}
\DeclareMathOperator*{\argmin}{arg\,min}

\usepackage{amsthm}
\newtheorem{theorem}{Theorem}
\newtheorem{proposition}{Proposition}

\newtheorem{corollary}{Corollary}

\newtheorem{lemma}{Lemma}[section]

\usepackage{dsfont}
\def\one{\mathds{1}}

\usepackage{amssymb}
\usepackage{makecell}
\usepackage{caption}
\usepackage{graphicx}

\usepackage{relsize}
\usepackage{subcaption}

\def\th{^\text{th}}
\def\T{{\mathsf T}}

\newcommand{\R}{\mathbb{R}}

\def\ball{{\mathsf B}}

\newcommand{\op}{{\rm op}}

\renewcommand{\P}{\mathbb{P}}
\newcommand{\E}{\mathbb{E}}

\def\de{\mathrm{d}}

\def\cS{{\mathcal S}}

\def\cN{{\mathcal N}}

\def\htheta{{\hat{\theta}}}
\def\hthetalambda{{\hat{\theta}_\lambda}}

\def\thetagdT{{\theta_T}}

\def\thetaegdT{{\theta_T}}
\newcommand{\thetamdT}{{\theta_T}}

\def\boldmu{{\mu}}

\def\risk{{\mathrm{Risk}}}

\def\sG{{\mathrm{sG}}}

\def\hSigma{{\widehat{\Sigma}}}
\def\trhSigma{{\mathrm{tr}(\widehat{\Sigma})}}
\def\FhSigma{{\|\widehat{\Sigma}\|_F}}
\def\ophSigma{{\|\widehat{\Sigma}\|_\op}}

\def\dist{{\mathrm{dist}}}
\def\proj{\Pi}

\def\cC{{\mathcal{C}}}
\def\cE{{\mathcal{E}}}

\newcommand{\bA}{{\mathsf{A}}}

\newcommand{\KL}{{\mathrm{KL}}}
\newcommand{\DKL}{{D_\KL}}

\newcommand{\Dphi}{{D_\phi}}

\def\simplex{{\Delta_d}}
\def\interior{{\mathrm{int}}}

\def\prox{{\mathrm{prox}}}

\def\hR{{\widehat R}}

\usepackage{ragged2e}
\usepackage{etoolbox}
\usepackage{easyReview}
\usepackage{cleveref}
\usepackage{lineno}
\usepackage[shortcuts]{extdash}

\title{Basic Inequalities for First-Order Optimization Algorithms with \\  
  Applications to Statistical Risk Analysis} 

\author{
  Seunghoon Paik\textsuperscript{\textdagger} \quad 
  Kangjie Zhou\textsuperscript{\textdaggerdbl} \quad 
  Matus Telgarsky\textsuperscript{\textsection} \quad 
  Ryan J. Tibshirani\textsuperscript{\textdagger} \\
  {\small
    \textsuperscript{\textdagger}University of California, Berkeley \quad
    \textsuperscript{\textdaggerdbl}Columbia University \quad
		\textsuperscript{\textsection}New York University
  }
}
\date{}

\begin{document}
\maketitle
\vspace{-20pt}

\begin{abstract}
In this work, we introduce \emph{basic inequalities} for first-order iterative
optimization algorithms, forming a simple yet versatile framework which connects
implicit and explicit regularization. Building on related comparison
inequalities for optimization iterates that already exist in the literature, we
extend and unify these arguments to produce a general framework, which can be
used as a tool for statistical analysis. In more detail, let $f$ denote the
objective function to be optimized. Given a first-order iterative algorithm
initialized at $\theta_0$, with current iterate $\theta_T$, the basic inequality
upper bounds $f(\theta_T) - f(z)$ for any reference point $z$ in terms of the
accumulated step sizes, and the distances between $\theta_0$, $\theta_T$, and
$z$. These distances are measured in a geometry inherent to the optimization
algorithm, which then translates into a notion of regularization being applied
across the path of iterates.

In addition to refining existing results on gradient descent, we provide new
results for mirror descent and other first-order methods.  We then show how to
use these basic inequalities to derive elementary yet useful bounds on the
prediction risk of early-stopped gradient descent and exponentiated gradient
descent iterates in generalized linear models. We also supplement these findings
with numerical experiments.
\end{abstract}

\section{Introduction}

This paper introduces \emph{basic inequalities} for first-order optimization
algorithms, which connects implicit and explicit notions of regularization. We
consider an optimization problem 
\[
\min_{\theta \in \R^d} \; f(\theta)
\]
for an objective function $f$, consider a first-order optimization method which
produces iterate $\theta_T$ at iteration $T$, and study the relationship between
$\theta_T$ and \smash{$\hthetalambda$}, where the latter solves 
\[
\min_{\theta \in \R^d} \; f(\theta) + \lambda g(\theta)
\]
for a regularizer $g$, and a tuning parameter $\lambda > 0$. By
construction, the solution \smash{$\hthetalambda$} is subject to traditional
regularization, explicitly shaped by $g$. In comparison, the iterate
\smash{$\theta_T$} is subject to implicit regularization, as we consider
iterative algorithms which enforce inductive bias along the optimization
trajectory, ultimately leading to specific structured solutions in
overparametrized settings. The study of iterate behavior at a finite iteration
$T$ is often called \emph{early stopping}.

The basic inequality derived in this paper upper bounds the gap
$f(\theta_T)-f(z)$ for any reference point $z$ in terms of the distances
between $\theta_0$, $\theta_T$, and $z$, and the accumulated step sizes, which
represent the total elapsed ``time'' in the optimization process. As an example, 
for gradient descent initialized at the origin, $\theta_0=0\in\R^d$, and a
constant step size $\eta>0$, we have the following basic inequality:    
\begin{equation}
\label{eq:basic_ineq_opt}
f(\theta_T) - f(z) \leq \frac{1}{2\eta T} \Big( \|z\|_2^2 - \|\theta_T - z\|_2^2
\Big),
\end{equation}
for any $z\in\R^d$ and $T \geq 1$. Meanwhile, for an explicitly regularized
solution, we have the following inequality: 
\begin{equation}
\label{eq:basic_ineq_reg}
f(\hthetalambda) - f(z) \leq \lambda \Big( g(z) - g(\hthetalambda) \Big),
\end{equation}
for any $z\in\R^d$ and $\lambda > 0$, obtained by rearranging the statement
that $f(\theta) + \lambda g(\theta)$ is minimized at \smash{$\hthetalambda$}.  

The inequality in \eqref{eq:basic_ineq_reg} is often used as a starting point in
the statistical analysis of regularized estimators. For example, consider a
regression loss \smash{$f(\theta) = \frac{1}{2n} \|Y-X\theta\|_2^2$} for a
response $Y\in\R^n$ and features $X\in\R^{n\times d}$, and $\ell_1$ penalty
$g(\theta) = \|\theta\|_1$; we can instantiate \eqref{eq:basic_ineq_reg} at the
regression coefficients $z=\theta^*$ in a population linear model for $Y|X$, and
rearrange further to yield  
\[
\frac{1}{2n} \| X\hthetalambda - X\theta^* \|_2^2\leq \frac{1}{n}
\langle X^\T \epsilon, \hthetalambda - \theta^* \rangle + \lambda 
\Big( \|\theta^*\|_1 - \|\hthetalambda\|_1\Big),
\]
where $\epsilon = Y-X\theta^* \in \R^n$ is a noise vector, and $\langle a, b
\rangle = a^\T b$. The above is referred to as the basic inequality in the 
high-dimensional statistics literature \citep{buhlmann2011statistics}, and it is
a central tool for analyzing the risk of \smash{$\hthetalambda$}, as an estimate
of $\theta^*$. Given the close resemblance between \eqref{eq:basic_ineq_opt} and
\eqref{eq:basic_ineq_reg}, one might imagine that the former can thus also be
used as a tool for analyzing the risk of the early-stopped estimate
$\theta_T$. This will be a large part of the focus in our paper. 

\subsection{Summary of contributions}

The inequality \eqref{eq:basic_ineq_opt} as well as numerous variations on this
idea can be found in the optimization literature on convergence analyses of
iterative algorithms \citep{nesterov2003introductory, nemirovski2009robust,
  reddi2018convergence}, and the machine learning literature on implicit
regularization and generalization \citep{ji2019implicit, ji2020gradient,
  wu2024large, wu2025benefits}. In this paper, we isolate and highlight a
certain technique for constructing such inequalities, and show that this
provides a framework for the statistical analysis of early-stopped estimates. A
summary of our specific contributions is as follows.   

\begin{itemize}
	\item \emph{A single analytical framework.}
	We introduce \emph{basic inequalities} for the iterates $\theta_T$ produced
  by popular first-order algorithms, including gradient descent and mirror
  descent (Section \ref{sec:basic-ineq}). These serve as
  well\-/rounded tools for understanding implicit regularization, as detailed in 
  below.   
	
\item \emph{Training envelopes and dynamics.} 
	Starting from the basic inequality, we derive a training envelope that
  provides lower and upper bounds on the combined loss and penalty of the
  iterates relative to explicitly regularized solutions (Section
  \ref{sec:basic-ineq}).  We also show how the basic inequality enables us to
  characterize the asymptotic behavior of solutions (as $T \to \infty$).
	
\item \emph{Statistical risk bounds.}  
	We use the basic inequality to derive high-probability bounds on the
  prediction risk of early-stopped iterates, and demonstrate that they match the
  rates obtained through traditional analysis of their explicit-regularized 
  counterparts. Specifically, given a sample size $n$, dimension $d$, and 
  regularity bound $b$, our framework yields the following results (where
  \smash{$\tilde O(\cdot)$} ignores log factors):
	\begin{itemize}
  \item For gradient descent on generalized linear model (GLM) losses, we obtain
    an excess risk bound of \smash{$\tilde O(b\sqrt{d/n})$}, matching the rate
    obtained by explicit ridge regularization (Section \ref{sec:glm}).    
		
  \item For exponentiated gradient descent on GLM losses, we instead obtain a
    bound of \smash{$\tilde O(\sqrt{(b\log d)/n})$}, matching the rate obtained
    by explicit KL regularization (Section \ref{sec:model-agg}).   
		
  \item For exponentiated gradient descent on linear losses in a randomized
    aggregation problem setting, we obtain an excess risk bound of
    \smash{$\tilde O(\sqrt{(b \log d)/n)}$}, and we note that this also applies
    to explicit KL-regularized solutions because here early stopping and
    explicit regularization actually result in identical paths (Section
    \ref{sec:alquier}).    
	\end{itemize} 
	
\item \emph{Experiments.}
	Complementing our theory, we run experiments with gradient descent and 
  exponentiated gradient descent on linear, logistic, and Poisson regression
  tasks. Across both underparameterized and overparameterized regimes, 
	the results show strong empirical alignment between implicit and explicit
  regularization in terms of training dynamics, prediction risk, and solution
  paths (Section \ref{sec:experiments}). 
\end{itemize}

Throughout this paper, we make use of standard definitions and notation, and we
provide an overview in the appendix for completeness. Furthermore, with the
exception of some key proofs for the basic inequalities, all proofs will be
deferred to the appendix.   

\subsection{Related work}

In what follows, as we present our results, we will then highlight and discuss the 
relevant literature in detail. Here we give a broader overview of work
related to our paper. At a high level, the term implicit regularization refers
to the phenomenon in which an optimization algorithm enforces stronger
properties (typically, enforces greater regularity) than what is suggested
purely by the perspective of loss minimization. A major interest in the
literature on implicit regularization has been the study of solutions an
asymptotic sense, i.e., as the given algorithm approaches a limit point. This
idea can be traced back to the literature on boosting and the study of
max-margin classifiers \citep{schapire1998boosting, rosset2004boosting, 
telgarsky2013margins}.
  
In the deep learning age, the observation that neural network models can and
will often generalize despite significant overparametrization has reinvigorated
interest in implicit regularization \citep{neyshabur2014search}. In the vein
of earlier boosting analyses, much progress has been made on understanding the
limiting behavior of gradient descent \citep{soudry2018implicit, ji2019implicit,
  ji2020gradient}, and mirror descent \citep{gunasekar2018characterizing,
  sun2023unified}.   

Distinct from the asymptotic behavior of an algorithm, early stopping considers
the behavior of iterates at a finite but variable number of iterations
$T$. Originally viewed as a heuristic to prevent overfitting in neural networks 
\citep{prechelt2002early}, a formal understanding began to develop, with
boosting again being the primary motivating framework early on
\citep{buhlmann2003boosting, zhang2005boosting}. There is by now a significant 
body of work connecting and comparing early-stopped gradient descent to explicit
$\ell_2$ (or ridge) regularization \citep{yao2007early, raskutti2014early,
  wei2017early, suggala2018connecting, ali2019continuous, ali2020implicit, 
  zou2021benefits, sonthalia2024regularization, wu2024large, wu2025benefits}. 

Similarly, a central focus in the current paper is to connect gradient descent
to ridge regularization, and more generally, mirror descent to Bregman
divergence regularization (Section \ref{sec:basic-ineq}). Our results are less
refined than what can be found in the literature, which analyzes the precise
behavior of descent iterates for particular problems---particular loss
functions, often paired with assumptions on the data at hand. That said, our
take is more general: we provide a general recipe for connecting optimization
iterates and explicit regularization through what we call \emph{basic
  inequalities}. This is based purely on optimization-theoretic arguments, and 
therefore only relies on high-level properties of the loss (e.g., Lipschitz
smoothness). The inequalities are deterministic, and while general and hence
somewhat coarse, they are tight enough to enable meaningful corollaries about
the statistical risk of early-stopped estimates in various settings (Sections
\ref{sec:glm}--\ref{sec:alquier}). Finally, the framework for basic inequalities
extends beyond gradient and mirror descent to other first-order algorithms
(Section \ref{sec:other-algos}), and we believe will be a fruitful platform for
future development.  

\section{Basic inequalities}
\label{sec:basic-ineq} 

This section introduces basic inequalities for gradient descent and mirror
descent, and shows how they can be used to understand training dynamics. While
the basic inequality for mirror descent strictly generalizes that for gradient
descent, we cover gradient descent first. This is done for the sake of
exposition, and also because the $\ell_2$ geometry underlying gradient descent
leads to some special corollaries.

\subsection{Gradient descent}

Given a differentiable loss function $f:\R^d\to \R$, gradient descent with an
initial point $\theta_0\in\R^d$, and step sizes
\smash{$\{\eta_t\}_{t=0}^\infty$}, generates iterates for $t=0,1,2,\dots$ via   
\begin{equation} 
\label{eq:gd-iterate}
\theta_{t+1} = \theta_t - \eta_t \nabla f(\theta_t).
\end{equation}
Theorem \ref{thm:gd-basicineq} presents our first result, a basic inequality
for gradient descent, which bounds the objective value of the iterate
$\theta_T\in\R^d$ with respect to any reference point $z\in\R^d$.   

\begin{theorem}
\label{thm:gd-basicineq}
Let $f:\R^d\to \R$ be convex, differentiable, and $L$-smooth, for
$L>0$. Consider gradient descent \eqref{eq:gd-iterate} with step sizes
$\eta_t\in (0,1/L]$, $t=0,1,2,\dots$. For any reference point $z\in\R^d$,     
and any iteration $T \geq 1$, 
\[
f(\theta_T) - f(z) \leq \frac{1}{2\sum_{t=0}^{T-1} \eta_t} \Big( \|\theta_0 -
z\|_2^2 - \|\theta_T - z\|_2^2 \Big). 
\]
In particular, for a constant step size $\eta_t = \eta$, $t=0,1,2,\dots$, this
simplifies to 
\begin{equation} 
\label{eq:gd-basicineq}
f(\theta_T) - f(z) \leq \frac{1}{2\eta T} \Big( \|\theta_0 - z\|_2^2 -
\|\theta_T - z\|_2^2 \Big). 
\end{equation}
\end{theorem}

\begin{proof}
The proof can be broken down into three steps.

\medskip\noindent
\underline{\smash{Step 1:}} 
\emph{Tracking the proximity difference across adjacent iterations.} 
Measuring proximity via the Euclidean distance, note that  
\[
\|\theta_t - z\|_2^2 - \|\theta_{t+1} - z\|_2^2 
= \|\theta_t - z\|_2^2 - \|\theta_t - \eta_t \nabla f(\theta_t) - z\|_2^2
= 2\eta_t \langle \nabla f(\theta_t), \theta_t - z\rangle - \eta_t^2 \|\nabla
f(\theta_t)\|_2^2. 
\] 
    
\smallskip\noindent
\underline{\smash{Step 2:}}
\emph{Bounding the objective difference $f(\theta_t) - f(z)$.} 
By convexity, we have $f(\theta_t) - f(z) \leq \langle \nabla f(\theta_t),
\theta_t - z\rangle$. Substituting this into the result from Step 1, 
\[
2\eta_t \Big(f(\theta_t) - \frac{\eta_t}{2} \|\nabla f(\theta_t) \|_2^2 - f(z) 
\Big) \leq \|\theta_t - z\|_2^2 - \|\theta_{t+1} - z\|_2^2. 
\]
By $L$-smoothness of $f$ and $\eta_t \in (0,1/L]$, the descent lemma (Lemma
\ref{lem:gd-descending}) guarantees   
\[
f(\theta_{t+1})
\leq f(\theta_t) - \eta_t\Big(1-\frac L2 \eta_t\Big) \|\nabla f(\theta_t)\|_2^2 
\leq f(\theta_t) - \frac{\eta_t}{2} \|\nabla f(\theta_t)\|_2^2.
\]
This ensures \smash{$ f(\theta_T) \leq f(\theta_t) - (\eta_t/2) \|\nabla 
  f(\theta_t)\|_2^2$} for $t<T$. Using this in the second-to-last display, 
\[
2\eta_t(f(\theta_T) - f(z)) \leq \|\theta_t - z\|_2^2 - \|\theta_{t+1} -
z\|_2^2. 
\]

\smallskip\noindent 
\underline{\smash{Step 3:}}
\emph{Aggregating bounds over iterations.} 
Summing the result of Step 2 over $t<T$ gives a telescoping sum: 
\[
2\sum_{t=0}^{T-1} \eta_t (f(\theta_T) - f(z)) \leq \|\theta_0 - z\|_2^2 -
\|\theta_T - z\|_2^2, 
\]
which completes the proof after rearrangement.
\end{proof}

Some readers will see \eqref{eq:gd-basicineq} as a familiar result, as it is
just a few steps away from a ``textbook''  result for gradient descent. To
derive the standard $1/T$ convergence rate on the suboptimality gap of 
gradient descent under Lipschitz smoothness, we can set $z$ to be a minimizer
of $f$, and then drop the negative term involving $\|\theta_T - z\|_2^2$ on 
the right-hand side. This gives the familiar bound \smash{$f(\theta_T) - f^\star
  \leq \frac{b^2}{2\eta T}$}, where $f^\star$ is the optimal objective value, and
$b$ bounds the distance between $\theta_0$ and the solution set. 

Inspired by use of similar techniques in \citet{ji2019implicit, ji2020gradient,
  wu2024large, wu2025benefits}, in \eqref{eq:gd-basicineq} we recognize the
importance of keeping $z\in\R^d$ as a free variable, and maintaining both 
initial $\|\theta_0 - z\|_2^2$ and current $\|\theta_T - z\|_2^2$ distances in
the bound. Though these are seemingly simple choices, the specific form of 
\eqref{eq:gd-basicineq} will lead to numerous insights on the behavior of
gradient descent iterates, as we will see in what follows. As our first example,
we show how the basic inequality \eqref{eq:gd-basicineq} leads to a
ridge-regularized training envelope for gradient descent iterates.      

\begin{corollary}
\label{cor:gd-envelope}
Under the assumptions of Theorem \ref{thm:gd-basicineq}, consider gradient
descent \eqref{eq:gd-iterate} with a constant step size $\eta_t = \eta \in
(0,1/L]$, $t=0,1,2,\dots$. Then for any $T \geq 1$ and $\lambda_T = 1/(\eta T)$,
\begin{equation} 
\label{eq:gd-envelope}
\min_{z\in\R^d} \Big( f(z) + \frac{\lambda_T}{4} \|\theta_0 - z\|_2^2 \Big)
\leq f(\theta_T) + \frac{\lambda_T}{4} \|\theta_0 - \theta_T\|_2^2
\leq \min_{z\in\R^d} \Big( f(z) + \lambda_T \|\theta_0 - z\|_2^2 \Big). 
\end{equation}
\end{corollary}

\begin{proof}
Observe that
\begin{align*}
\|\theta_0 - z\|_2^2 - \|\theta_T - z\|_2^2
&= 2\langle\theta_0-z, \theta_0-\theta_T\rangle - \|\theta_0-\theta_T\|_2^2 \\
&\leq 2 \|\theta_0-z\|_2 \|\theta_0-\theta_T\|_2 - \|\theta_0-\theta_T\|_2^2 \\
&\leq 2\|\theta_0-z\|_2^2 - \frac{1}{2}\|\theta_T-\theta_0\|_2^2,
\end{align*}
where the last line uses Young's inequality, $2ab \leq ca^2 + b^2/c$, with
$c=2$. Plugging this into \eqref{eq:gd-basicineq}, and then minimizing over 
$z\in\R^d$, gives the upper bound in \eqref{eq:gd-envelope}. The lower bound is
immediate. 
\end{proof}

The result in \eqref{eq:gd-envelope} translates the ``elapsed time'' $\eta T$ in
gradient descent to an effective ridge regularization parameter $\lambda_T$.  
For simplicity, set $\theta_0=0$, and let \smash{$\htheta_\lambda \in \R^d$}
be the (unique) minimizer of the ridge-penalized criterion $f(z) + \lambda
\|z\|_2^2$ over $z\in\R^d$. Rephrased, the result in \eqref{eq:gd-envelope} says
that for all $T \geq 1$,  
\[
f(\htheta_{\lambda_T/4}) + \frac{\lambda_T}{4} \|\htheta_{\lambda_T/4}\|_2^2 
\leq f(\theta_T) + \frac{\lambda_T}{4} \|\theta_T\|_2^2 
\leq f(\htheta_{\lambda_T}) + \lambda_T \|\htheta_{\lambda_T}\|_2^2.  
\]
In other words, while (by definition) the iterate $\theta_T$ must be suboptimal
with respect to the ridge-penalized criterion at a penalty factor $\lambda_T/4$, 
its achieved criterion value is no more than the optimal ridge-penalized
criterion at an inflated penalty factor $\lambda_T$. These lower and upper
bounds together provide an envelope for the ridge-penalized criteria achieved by
gradient descent across the full path, providing interesting evidence of the  
implicit $\ell_2$ regularity present in the gradient descent iterates. Notably,
this is true for any convex and Lipschitz smooth loss $f$.  

The next result demonstrates how the basic inequality is powerful enough to 
recover various standard facts about the training dynamics of gradient
descent. Our approach and presentation here is inspired by the work of
\citet{lemaire1996asymptotical} on gradient flow.   

\begin{corollary}
\label{cor:gd-optlemaire} 
Under the assumptions of Theorem \ref{thm:gd-basicineq}, consider gradient
descent \eqref{eq:gd-iterate} with step sizes \smash{$\eta_t \in (0,1/L]$},
$t=0,1,2,\dots$. The following holds, where \smash{$\proj_S(u) = \argmin_{s\in
    S} \|u-s\|_2$} denotes the projection of a point $u$ onto a set $S$, and 
\smash{$\dist_S(u) = \min_{s\in S} \|u-s\|_2$} denotes the distance from $u$ to
$S$. 

\begin{enumerate}
\item (Objective convergence.)
  Let \smash{$f^\star = \inf_{\theta\in\R^d} f(\theta) \in [-\infty,\infty)$}.
  If \smash{$\sum_{t=0}^\infty \eta_t = \infty$}, then \smash{$\lim_{t\to\infty}
    f(\theta_t) = f^\star$}.
		
\item (Nonincreasing distance to solution set.)
  Define the solution set \smash{$S =\{\theta \in\R^d: f(\theta) = f^\star\}$}.
  Notice that $S$ is closed and convex. If $S\neq\emptyset$, then    
  \[
  \{\|\theta_t-s\|_2\}_{t=0}^\infty \; \text{is nonincreasing for any $s\in S$},
  \quad\text{and thus},\quad 
  \{\dist_S(\theta_t)\}_{t=0}^\infty\;\text{is nonincreasing}.
  \]

\item (Iterate convergence.)
  If $S\neq\emptyset$ and \smash{$\sum_{t=0}^\infty \eta_t = \infty$}, then 
  \smash{$\lim_{t\to\infty} \theta_t = \theta_\infty \in S$}. Furthermore, 
  \[ 
  \|\theta_\infty-\proj_S(\theta_0)\|_2 \leq \dist_S(\theta_0), 
  \quad\text{and thus},\quad 
  \|\theta_\infty-\theta_0\|_2 \leq 2 \cdot \dist_S(\theta_0).
  \] 
		
\item (Projection onto the solution set.)
   If $S\neq\emptyset$, \smash{$\sum_{t=0}^\infty \eta_t = \infty$}, and $S$ is
   affine, then \smash{$\theta_\infty = \proj_S(\theta_0)$}.  
\end{enumerate}
\end{corollary}

\begin{proof}
We prove each part separately.

\medskip\noindent
\underline{Part 1.} 
By the first result in Theorem \ref{thm:gd-basicineq}, we know that 
\smash{$f(\theta_T) \leq f(z) + \|\theta_0 - z\|_2^2/(2\sum_{t=0}^{T-1}
  \eta_t)$} for any $z \in \R^d$. Since \smash{$\sum_{t=0}^\infty \eta_t =
  \infty$}, we have \smash{$\limsup_{T\to\infty} f(\theta_T) \leq f(z)$}, and
taking an infimum over $z$ gives the result.
	
\medskip\noindent
\underline{Part 2.}
By the first result in Theorem \ref{thm:gd-basicineq}, for any $s \in S$, 
\[
\|\theta_T - s\|_2^2 \leq \|\theta_0 - s\|_2^2 + 2 \bigg(\sum_{t=0}^{T-1}  
\eta_t\bigg) ( f(s) - f(\theta_T)) \leq \|\theta_0 - s\|_2^2,
\]
where the last inequality is due to $f(s) \leq f(\theta_T)$. Now, $\theta_T$ can
also be seen as the result of running gradient descent from an initial point
$\theta_t$ for $T-t$ iterations, and hence by the same argument $\|\theta_T -
s\|_2 \leq \|\theta_t - s\|_2$, which proves the first claim. Taking an infimum
over $s \in S$ proves the second claim.
	
\medskip\noindent
\underline{Part 3.}
Since \smash{$\{\|\theta_t-s\|_2\}_{t=0}^\infty$} is nonincreasing and bounded
below by zero, the Bolzano-Weierstrass theorem implies that there exists a
subsequence which converges, i.e., \smash{$\lim_{k\to\infty} \theta_{t_k} =
  \theta_\infty$}. By part 1 and continuity of $f$, we know that $\theta_\infty
\in S$. Then once again as \smash{$\{\|\theta_t-s\|_2\}_{t=0}^\infty$} is
nonincreasing, we know that \smash{$\lim_{t\to\infty} \theta_t =
  \theta_\infty$}, the first claim. The second claim follows from
\smash{$\{\|\theta_t-s\|_2\}_{t=0}^\infty$} being nonincreasing for $s = 
\proj_S(\theta_0)$. 

\medskip\noindent
\underline{Part 4.}
Let $p = \proj_S(\theta_0)$ and $v = p - \theta_\infty$. Assume $v \neq 0$. For
any $c \geq 0$, define $\beta_c = p + c \cdot \dist_S(\theta_0) \cdot
v/\|v\|_2$. As $S$ is affine, we know that $\beta_c \in S$. Then by parts 2 and
3, we know that $\|\theta_\infty - \beta_c\|_2 \leq \|\theta_0 -
\beta_c\|_2$. Since the three points $\beta_c, \theta_\infty, p$ are collinear
by construction,  
\[
\|\theta_\infty - \beta_c\|_2
= \|\theta_\infty - p\|_2 + \|p - \beta_c\|_2
= \|v\|_2 + c \cdot \dist_S(\theta_0).
\]
Meanwhile, by the Pythagorean theorem,
\[
\|\theta_0 - \beta_c\|_2^2
= \|\theta_0 - p\|_2^2 + \|p - \beta_c\|_2^2
= \dist_S(\theta_0)^2 + (c \cdot \dist_S(\theta_0))^2
= (1+c^2) \cdot \dist_S(\theta_0)^2.
\]
Therefore, we finally have
\[
\|v\|_2 + c \cdot \dist_S(\theta_0) 
= \|\theta_\infty - \beta_c\|_2
\leq \|\theta_0 - \beta_c\|_2
= \sqrt{1+c^2} \cdot \dist_S(\theta_0).
\]
Thus \smash{$\|v\|_2 \leq (\sqrt{1+c^2} - c) \cdot \dist_S(\theta_0)$}. Taking
$c \to \infty$ implies $\|v\|_2=0$, which means $p = \theta_\infty$. 
\end{proof}

We note that part 4 of the above corollary reproduces a well-known ``folklore'' 
result in overparameterized linear regression: gradient descent initialized at
the origin converges to the minimum $\ell_2$ norm solution. More generally, it
converges to the projection of the initial point onto the solution set. This
also covers any loss $f$ that depends on $\theta$ only through a linear
transformation $X\theta$, such as generalized linear model losses (as introduced
formally in Section \ref{sec:glm}).      

\subsection{Mirror descent} 

Mirror descent is a generalization of gradient descent where Euclidean distance
is replaced by the \emph{Bregman divergence} induced by a convex, differentiable
function $\phi:\Omega\to\R$. For points $u,v \in \interior(\Omega)$
(the interior of the domain $\Omega$, or relative interior in the case that $\Omega$ is not full-dimensional)
this is defined by
\[ 
\Dphi(u,v) = \phi(u) - \phi(v) - \langle \nabla\phi(v), u-v \rangle.
\]
For a constraint $\cC \subseteq \Omega$, an initial point $\theta_0 \in \cC \cap
\interior(\Omega)$, and step sizes \smash{$\{\eta_t\}_{t=0}^\infty$}, the
mirror descent update is 
\begin{equation} 
\label{eq:md-iterate}
\theta_{t+1} = \argmin_{\theta\in\cC} \Big( \eta_t \langle \nabla f(\theta_t),
\theta\rangle + \Dphi(\theta, \theta_t) \Big). 
\end{equation}
This has an equivalent two-stage formulation, where we first perform a gradient
step in a transformed space given by the ``mirror map'' $\nabla\phi$, and then
take a Bregman projection onto $\cC$:
\[
\nabla\phi(\tilde\theta_{t+1}) = \nabla\phi(\theta_t) - \eta_t \nabla
f(\theta_t), \qquad  
\theta_{t+1} = \argmin_{\theta\in\cC} \, \Dphi(\theta, \tilde\theta_{t+1}).    
\]

A popular example of mirror descent is the exponentiated gradient descent
algorithm, which is designed for optimization over the probability simplex,
\smash{$\Delta_d=\{a\in\R^d: a_i\geq 0, \, \sum_{i=1}^d a_i = 1\}$}.
By choosing the negative entropy function \smash{$\phi(u) = \sum_{i=1}^d u_i
  \log u_i$}, which is 1-strongly convex with respect to $\|\cdot\|_1$, on
$\Omega=\cC=\Delta_d$, the induced Bregman divergence becomes 
the Kullback-Leibler (KL) divergence, \smash{$\DKL(a,b)=\sum_{i=1}^d
  a_i\log(a_i/b_i)$}. The corresponding mirror descent update, where $\odot$ 
denotes coordinatewise multiplication, is
\begin{equation}
\label{eq:egd-iterate}
\tilde\theta_{t+1} = \theta_t \odot \exp(-\eta_t\nabla f(\theta_t)),
\qquad
\theta_{t+1}=\tilde\theta_{t+1} / \|\tilde\theta_{t+1}\|_1.
\end{equation}

Mirror descent can also recover projected gradient descent with constraint set
$\cC$, which of course reduces to ordinary gradient descent when $\cC=\R^d$. 
Taking \smash{$\phi(u) = \frac12\|u\|_2^2$}, the induced Bregman divergence is
simply \smash{$\Dphi(u,v)=\frac12\|u-v\|_2^2$}, and the mirror descent update is 
\begin{equation} 
\label{eq:pgd-iterate}
\theta_{t+1} = \proj_{\cC}(\theta_t - \eta_t \nabla f(\theta_t)),
\end{equation}
where recall $\proj_{\cC}$ denotes Euclidean projection onto $\cC$.

Below we introduce a basic inequality for mirror descent iterates in
\eqref{eq:md-iterate}. Just as mirror descent generalizes gradient descent, 
Theorem \ref{thm:md-basicineq} generalizes Theorem \ref{thm:gd-basicineq}. 

\begin{theorem}
\label{thm:md-basicineq}
Let $\Omega,\cC$ be closed, convex sets in $\R^d$ with nonempty interior, and 
$\cC\subseteq\Omega$. Assume $f:\Omega\to\R$ is convex on $\cC$, and
$L$-smooth with respect to a norm $\|\cdot\|$ on
$\cC\cap\interior(\Omega)$, for $L>0$. Assume $\phi:\Omega\to\R$ is of
Legendre type, and $\alpha$-strongly convex with respect to $\|\cdot\|$ on
$\Omega$, for $\alpha>0$. Consider mirror descent \eqref{eq:md-iterate},
initialized at $\theta_0\in\cC\cap\interior(\Omega)$, with step sizes
$\eta_t\in(0,\alpha/L]$, $t=0,1,2,\dots$. For any reference point $z\in\cC$,
and any $T \geq 1$,  
\[
f(\theta_T) - f(z) \leq \frac{1}{\sum_{t=0}^{T-1}\eta_t}\Big(\Dphi(z, \theta_0)
- \Dphi(z, \theta_T)\Big). 
\]
In particular, for a constant step size $\eta_t = \eta$, $t=0,1,2,\dots$, this
simplifies to 
\begin{equation}
\label{eq:md-basicineq} 
f(\theta_T) - f(z) \leq \frac{1}{\eta T} \Big(\Dphi(z, \theta_0) - \Dphi(z,
\theta_T)\Big). 
\end{equation}
\end{theorem}

\begin{proof}
The proof parallels that from the gradient descent case.

\medskip\noindent
\underline{\smash{Step 1:}} 
\emph{Tracking the proximity difference across adjacent iterations.} 
Measuring proximity via the Bregman divergence, the three-point identity for the
Bregman divergence states that  
\[ 
\Dphi(z, \theta_{t+1}) + \Dphi(\theta_{t+1}, \theta_t) - \Dphi(z, \theta_t) =
\langle \nabla\phi(\theta_t) - \nabla\phi(\theta_{t+1}), z -
\theta_{t+1}\rangle. 
\]
Note $\langle \eta_t \nabla f(\theta_t) + \nabla\phi(\theta_{t+1}) -
\nabla\phi(\theta_t), z - \theta_{t+1} \rangle \geq 0$, by the first-order
optimality condition for $\theta_{t+1}$ in \eqref{eq:md-iterate}. Rearranging
terms, this implies $\eta_t\langle \nabla f(\theta_t), \theta_{t+1}-z\rangle
\leq \langle \nabla\phi(\theta_{t+1})-\nabla\phi(\theta_t),
z-\theta_{t+1}\rangle$, and combining with the above identity,    
\[
\eta_t \langle \nabla f(\theta_t), \theta_{t+1} - z\rangle
\leq \Dphi(z,\theta_t) - \Dphi(z, \theta_{t+1}) - \Dphi(\theta_{t+1}, \theta_t).
\]

\smallskip\noindent
\underline{\smash{Step 2:}}
\emph{Bounding the objective difference $f(\theta_t) - f(z)$.} 
By convexity, 
\[
f(\theta_t) - f(z) \leq \langle\nabla f(\theta_t), \theta_t - z \rangle
= \langle\nabla f(\theta_t), \theta_{t}- \theta_{t+1} \rangle
+ \langle\nabla f(\theta_t), \theta_{t+1} - z \rangle.
\]
Combining this with the result of Step 1, we have
\[
\eta_t(f(\theta_t) - f(z))
\leq \eta_t \langle\nabla f(\theta_t), \theta_{t}- \theta_{t+1} \rangle
+ \Dphi(z,\theta_t) - \Dphi(z, \theta_{t+1}) - \Dphi(\theta_{t+1}, \theta_t).
\]
Meanwhile, the $L$-smoothness of $f$ gives
$f(\theta_{t+1}) \leq f(\theta_t) + \langle \nabla f(\theta_t),
\theta_{t+1} - \theta_t\rangle + (L/2) \|\theta_{t+1}-\theta_t\|^2$,
and the $\alpha$-strong convexity of $\phi$ implies $\Dphi(\theta_{t+1}, 
\theta_t) \geq (\alpha/2) \|\theta_{t+1}-\theta_t\|^2$. Substituting these into
the above display,  
\begin{align*}
\eta_t (f(\theta_{t+1}) - f(z))
&\leq \Dphi(z,\theta_t) - \Dphi(z, \theta_{t+1}) - \Big(\frac\alpha2 - \frac L2 
\eta_t \Big) \|\theta_{t+1}-\theta_t\|^2 \\
&\leq \Dphi(z,\theta_t) - \Dphi(z, \theta_{t+1}),
\end{align*}
where we use $\eta_t\leq \alpha/L$ in the last inequality. The descent lemma for
mirror descent (Lemma \ref{lem:md-descending}) shows that $f(\theta_t)$ is
nonincreasing in $t$, thus we have 
\[
\eta_t (f(\theta_T) - f(z)) \leq \Dphi(z,\theta_t) - \Dphi(z, \theta_{t+1}). 
\]
    
\smallskip\noindent 
\underline{\smash{Step 3:}}
\emph{Aggregating bounds over iterations.} 
Summing the result of Step 2 over $t<T$ gives a telescoping sum: 
\[
\sum_{t=0}^{T-1} \eta_t (f(\theta_T) -f(z)) \leq \Dphi(z, \theta_0) -
\Dphi(z, \theta_T), 
\]
which completes the proof after rearrangement.
\end{proof}

Next we present the mirror descent analog of Corollary \ref{cor:gd-envelope}, on
a regularized training envelope. 

\begin{corollary}
\label{cor:md-envelope}
Under the assumptions of Theorem \ref{thm:md-basicineq}, consider mirror
descent \eqref{eq:md-iterate} with a constant step size $\eta_t = \eta \in
(0,\alpha/L]$, $t=0,1,2,\dots$. Further, assume that there exists $G>0$ such
that \smash{$\Dphi(z,\theta_0) \leq \frac{G}{2}\|\theta_0-z\|^2$} for any
$z\in\cC$ (note that this is implied by $\phi$ being $G$-smooth with respect to
$\|\cdot\|$ on $\cC$). Then for any $T \geq 1$ and $\lambda_T = 1/(\eta T)$,  
\begin{equation} 
\label{eq:md-envelope}
\min_{z\in\cC} \Big( f(z) + \frac{\alpha \lambda_T}{4} \|\theta_0 - z\|^2 \Big)
\leq f(\theta_T) + \frac{\alpha \lambda_T}{4} \|\theta_0 - \theta_T\|^2 
\leq \min_{z\in\cC} \Big( f(z) + \frac{(G+\alpha) \lambda_T}{2} \|\theta_0 -
z\|^2 \Big).    
\end{equation}
\end{corollary}

\begin{proof}
Observe that
\begin{align*}
\Dphi(z,\theta_0) - \Dphi(z,\theta_T) 
&\leq \frac{G}{2} \|\theta_0-z\|^2 - \frac{\alpha}{2} \|\theta_T-z\|^2 \\
&\leq \frac{G}{2} \|\theta_0-z\|^2 - \frac{\alpha}{2} \Big( \frac{1}{2}
\|\theta_T-\theta_0\|^2 - \|\theta_0-z\|^2 \Big) \\
&= \frac{G+\alpha}{2} \|\theta_0-z\|^2 - \frac{\alpha}{4}
\|\theta_T-\theta_0\|^2, 
\end{align*}
where the first line uses the given assumption on $\phi$ combined with
$\alpha$-strong convexity, and the second uses \smash{$\frac{1}{2}
  \|\theta_T-\theta_0\|^2 \leq \|\theta_T-z\|^2 + \|\theta_0-z\|^2$}, based on    
the triangle inequality and Young's inequality $2ab \leq a^2 + b^2$. We can plug 
this into \eqref{eq:md-basicineq}, and minimize over $z\in\cC$, to derive the  
upper bound in \eqref{eq:md-envelope}. As before, the lower bound is immediate.      
\end{proof}

We remark that \eqref{eq:md-envelope} recovers the previous gradient descent
result \eqref{eq:gd-envelope} when we take \smash{$\phi=\frac12\|\cdot\|_2^2$},
which implies $G=\alpha=1$. In this case, we have a tight envelope, with a 1:4
ratio of the penalty factors used in the comparison between the regularized
criteria achieved by gradient descent (penalty factor $\lambda_T/4$) and an
optimal solution (penalty factor $\lambda_T$).  In general, the envelope in
\eqref{eq:md-envelope} admits a ratio of 1:$(2G/\alpha+2)$, which can be loose
when $G$ is large relative to $\alpha$. In particular, for exponentiated
gradient descent \eqref{eq:egd-iterate} the function $\phi$ is the negative
entropy, in which case we have smoothness and strong convexity parameters $G=d$
and $\alpha=1$ with respect to the $\ell_1$ norm. This leads to a ratio of
1:$(2d+2)$, and therefore the upper bound in \eqref{eq:md-envelope} quickly
becomes uninformative in high dimensions. Our empirical results in Section
\ref{sec:experiments} suggest that this is overly conservative in
practice. Furthermore, using a reverse Pinsker-type inequality for the KL
divergence \citep{sason2015reverse}, we can derive an alternative upper bound
that (in certain regimes) has a better dimension dependence. This is presented 
below, and proved in the appendix.

\begin{corollary}
\label{cor:egd-envelope}
Under the assumptions of Theorem \ref{thm:md-basicineq}, consider now
exponentiated gradient descent \eqref{eq:egd-iterate} with a constant step size
$\eta_t = \eta \in (0,1/L]$, $t=0,1,2,\dots$ and initial point $\theta_0 = \pi$,
where $\pi = (1/d, \dots, 1/d) \in \Delta_d$ is the uniform distribution. Then
for any $T \geq 1$ and $\lambda_T = 1/(\eta T)$,  
\begin{multline} 
\label{eq:egd-envelope}
\min_{z\in\Delta_d} \Big( f(z) + \frac{\lambda_T}{4} \|\pi - z\|_1^2 \Big) 
\leq f(\theta_T) + \frac{\lambda_T}{4} \|\pi - \theta_T\|_1^2 \\
\leq \min_{z\in\Delta_d} \Big( f(z) + \lambda_T \cdot \min \Big\{ \frac{d+1}{2}
\|\pi - z\|_1^2, \, \frac{1}{2} \|\pi - z\|_1^2 + \frac{\log d}{2(1-1/d)} \|\pi
- z\|_1 \Big\} \Big).      
\end{multline}
\end{corollary}

Under mild conditions, we can reformulate the bound in this corollary using a
pure $\ell_1$ penalty. Provided we run enough iterations of exponentiated
gradient descent so that $\theta_T$ is bounded away from the uniform
distribution $\pi$, i.e., $\| \theta_T - \pi\|_1 \geq c$ for any constant $c>0$,
we can use this lower bound in \eqref{eq:egd-envelope} and the simple upper
bound $\|x-y\|_1 \leq 2$ for any $x,y \in \Delta_d$, to obtain
\begin{equation} 
\label{eq:egd-envelope2}
\min_{z\in\Delta_d} \Big( f(z) + \frac{c \lambda_T}{4} \|\pi - z\|_1 \Big)  
\leq f(\theta_T) + \frac{c \lambda_T}{4} \|\pi - \theta_T\|_1 
\leq \min_{z\in\Delta_d} \Big( f(z) + \lambda_T \Big[\frac{\log d}{2(1-1/d)} +
1\Big] \|\pi - z\|_1 \Big).     
\end{equation}
The envelope in \eqref{eq:egd-envelope2} clearly has a much better dimension 
dependence, with a ratio of 1:$O(\log{d})$ between the penalty factors for
exponentiated gradient descent and an optimal ($\ell_1$-regularized) solution.            

Lastly we present the mirror descent analog of Corollary
\ref{cor:gd-optlemaire}, on training dynamics. The proof is overall similar to
that of Corollary \ref{cor:gd-optlemaire}, and deferred to the appendix. 

\begin{corollary}
\label{cor:md-optlemaire}
Under the assumptions of Theorem \ref{thm:md-basicineq}, consider mirror descent 
\eqref{eq:md-iterate} with step sizes \smash{$\eta_t\in (0,\alpha/L]$},
$t=0,1,2,\dots$. The following holds, where \smash{$\proj^\phi_S(u) =
  \argmin_{s\in S} \Dphi(s,u)$} denotes the Bregman projection of a point $u$
onto a set $S$, and \smash{$\dist^\phi_S(u) = \min_{s\in S} \Dphi(s,u)$}
denotes the Bregman distance from $u$ to $S$.  

\begin{enumerate}
\item (Objective convergence.)
  Let \smash{$f^\star = \inf_{\theta\in\cC} f(\theta) \in [-\infty,\infty)$}.
  If \smash{$\sum_{t=0}^\infty \eta_t = \infty$}, then \smash{$\lim_{t\to\infty}
    f(\theta_t) = f^\star$}.
		
\item (Nonincreasing distance to solution set.)
  Define the solution set \smash{$S =\{\theta \in\cC: f(\theta) = f^\star\}$}.
  Notice that $S$ is closed and convex. If $S\neq\emptyset$, then    
  \[
  \{\Dphi(s,\theta_t)\}_{t=0}^\infty \; \text{is nonincreasing for any $s\in
    S$}, \quad\text{and thus},\quad 
  \{\dist^\phi_S(\theta_t)\}_{t=0}^\infty\;\text{is nonincreasing}.
  \]

\item (Iterate convergence.)
  Suppose $S\neq\emptyset$ and \smash{$\sum_{t=0}^\infty \eta_t = \infty$}. If
  either $S \cap \interior(\Omega) \neq \emptyset$, or $\Dphi$ is continuous in
  its second argument relative to $\interior(\Omega)$\footnote{To be precise,
    this means that for every sequence \smash{$\{y_k\}_{k=1}^\infty \subseteq
      \interior(\Omega)$} such that $y_k \to y$, we have $\Dphi(y,y_k) \to
    0$. Note this is not generally true for $\phi$ of Legendre type, see, e.g.,
    Remark 3.4 and Example 7.32  in \citet{bauschke1997legendre}.}, then 
  \smash{$\lim_{t\to\infty} \theta_t = \theta_\infty \in S$}. 
		
\item (Minimum Bregman divergence solution.)
   If $S\neq\emptyset$, \smash{$\sum_{t=0}^\infty \eta_t = \infty$}, and $S$ is
   affine with $S \subseteq \interior(\Omega)$, then \smash{$\theta_\infty =
     \proj^\phi_S(\theta_0)$}.   
\end{enumerate}
\end{corollary}

We note that part 4 of the above corollary reproduces a result of
\citet{gunasekar2018characterizing}; this covers losses $f$ which depend on 
$\theta$ only through a linear transformation $X\theta$, such as generalized
linear model losses, which we study in detail next. 

\section{Generalized linear models and ridge regularity}
\label{sec:glm} 

We now shift to a statistical perspective to analyze the prediction risk of the
iterates. This section focuses on generalized linear models (GLMs) and gives
a comparative analysis of two regularization schemes: implicit via gradient
descent and explicit via ridge regularization.  

We adopt a standard GLM setup with an identity sufficient statistic, where $Y
\in \R^n$ denotes the response vector and $X \in \R^{n\times d}$ the feature 
matrix. The negative log-likelihood function is thus (after rescaling by
\smash{$\frac1n$}):
\begin{equation} 
\label{eq:glm-loss}
f(\theta) = \frac1n \Big(-Y^\top X\theta + \bA(X\theta)\Big).
\end{equation}
for a map \smash{$\bA : \R^n \to \R$} which acts coordinatewise, as in
\smash{$\bA(v) = \sum_{i=1}^n A(v_i)$}. Here \smash{$A : \R \to \R$} is the   
cumulant generating function of the univariate exponential family associated
with the given GLM. Key examples are $A(u)=u^2/2$ for the Gaussian distribution 
(linear regression), $A(u) = \log(1+e^u)$ for the Bernoulli distribution
(logistic regression), and $A(u) = e^u$ for the Poisson distribution (Poisson
regression). An important property is that, for any exponential family, the
cumulant generating function $A$ is convex.

In what follows, we consider a fixed-X analysis of prediction risk in GLMs: we 
generally treat $X \in \R^{n\times d}$ as fixed (nonrandom). We denote by $P_i$
the distribution of each response $Y_i$, and write 
\[
\mu_i = \E[Y_i], \quad\text{and}\quad \epsilon_i = Y_i-\mu_i. 
\]
Unless otherwise stated, we assume that $Y_1,\dots,Y_n$ are independent
conditional on $X$.
We allow for model misspecification in the sense that we do not require the
mean $\mu_i$ to follow the canonical link of the given GLM. That is, when
analyzing linear regression we do not require the mean $\mu_i$ to be linear in
$x_i$ (the $i\th$ row of $X$), when analyzing logistic regression we do not
require the log odds $\log(\mu_i / (1-\mu_i))$ to be linear in $x_i$, and so on.  
We define the \emph{prediction risk} (or simply risk) of an estimator
\smash{$\htheta = \htheta(X,Y)$} as the expected negative log-likelihood
evaluated at a test copy \smash{$Y^* \in \R^n$} of the training response
vector, where each \smash{$Y^*_i \sim P_i$} (independent of $Y$):  
\begin{equation} 
\label{eq:glm-risk}
\risk(\htheta) = \frac1n \E \Big[ -(Y^*)^\top X\htheta + \bA(X\htheta)
\,\big|\, Y \Big]. 
\end{equation}
An important observation is that \smash{$\risk(\htheta) = \frac1n(-\mu^\top 
  X\htheta + \bA(X\htheta))$}, and we can therefore add and subtract $\epsilon$
to $\mu$ and expand to obtain: 
\[
\risk(\htheta) = f(\htheta) + \frac1n \epsilon^\top X \htheta.
\]
This expresses the prediction risk as the sum of the training error and a
stochastic term linear in the noise. This structure will allow us to apply the
basic inequality to derive high-probability risk bounds. 

\subsection{Risk analysis: ridge-regularized solution} 
\label{sec:risk-ridgeglm}

The ridge-regularized GLM estimator is defined by augmenting the GLM loss
\eqref{eq:glm-loss} with a squared $\ell_2$ penalty, denoted
\begin{equation} 
\label{eq:glm-ridge}
\hthetalambda = \argmin_{\theta\in\R^d} \Big( f(\theta) + \lambda
\|\theta\|_2^2 \Big).  
\end{equation}
for a regularization parameter $\lambda > 0$. (The solution exists and is unique
as the objective is strongly convex.) The purpose of this subsection is to
obtain relatively straightforward but nonetheless meaningful risk bounds for 
ridge-regularized GLM estimates, in order to compare what is achievable for
gradient descent estimates from the basic inequality, in the next subsection. 

As we will see, the risk bounds for explicitly- and implicitly-regularized
estimates will be comparable, and further, the strategies we use to derive
them will be similar. As a jumping off point for the ridge analysis, by
definition of \smash{$\hthetalambda$} in \eqref{eq:glm-ridge}, we know that for
any $\theta \in \R^d$,    
\[
f(\hthetalambda) + \lambda \|\hthetalambda\|_2^2 \leq f(\theta) +
\lambda\|\theta\|_2^2.
\]
This can be rewritten as
\[
f(\hthetalambda) - f(\theta) \leq \lambda \Big(\|\theta\|_2^2 -
\|\hthetalambda\|_2^2 \Big). 
\]
Adding \smash{$\frac1n \epsilon^\top X (\hthetalambda-\theta)$} to each side,
and recalling the decomposition of risk in \eqref{eq:glm-risk}, we get   
\begin{equation}
\label{eq:glm-risk-diff}
\risk(\hthetalambda) - \risk(\theta) \leq \frac1n \epsilon^\top X
(\hthetalambda-\theta) + \lambda \Big(\|\theta\|_2^2 - \|\hthetalambda\|_2^2
\Big).  
\end{equation}  
The next proposition develops this further to provide more salient deterministic
and high-probability bounds on the excess risk. Here and henceforth, a random
variable $Z$ is said to be sub-Gaussian with parameter \smash{$\sigma^2$},
written \smash{$Z\sim\sG(\sigma^2)$}, if \smash{$\E[\exp(t(Z-\E[Z]))] \leq
  \exp(\sigma^2 t^2/2)$} for all \smash{$t\in\R$}. 

\begin{proposition}
\label{prop:risk-ridgeglm}
For any \smash{$\lambda>0$} and any reference point \smash{$\theta\in\R^d$}, the 
prediction risk of \smash{$\hthetalambda$} in \eqref{eq:glm-ridge} satisfies 
\[
\risk(\hthetalambda) - \risk(\theta) \leq \frac{1}{2\lambda}
\Big\|\frac{X^\top\epsilon}{n}\Big\|_2^2 + 2\lambda\|\theta\|_2^2. 
\]
Furthermore, assume that $\epsilon_i$, $i=1,\dots,n$, are independent, each
$\epsilon_i \sim \sG(\sigma_i^2)$. Let \smash{$\sigma^2 =
\max_{i=1,\dots,n} \sigma_i^2$} and 
\smash{$\hSigma = \frac1n X^\top X$}. Then, for any $\delta>0$ and $b>0$, 
choosing 
\begin{equation} 
\label{eq:glm-ridge-lambda}
\lambda = \frac{\sigma}{2b\sqrt{n}} \sqrt{\trhSigma + 2\sqrt\delta\FhSigma + 
  2\delta\ophSigma},
\end{equation}
the excess risk of \smash{$\hthetalambda$} with respect the class of parameters 
with $\ell_2$ norm at most $b$ satisfies
\begin{equation} 
\label{eq:glm-ridge-riskhighprob}
\risk (\hthetalambda) - \min_{\theta: \|\theta\|_2\leq b} \risk (\theta)
\leq \frac{2b\sigma }{\sqrt{n}} \sqrt{\trhSigma + 2\sqrt\delta\FhSigma + 
  2\delta\ophSigma},
\end{equation}
with probability at least \smash{$1 - e^{-\delta}$}.
\end{proposition}

The bound \eqref{eq:glm-ridge-riskhighprob} is governed by the spectral
properties of the empirical covariance \smash{$\hSigma$}. Under the scaling 
\begin{equation} 
\label{eq:hsigma-spectral}
\trhSigma= O(d), \quad \FhSigma = O(\sqrt{d}), \quad \text{and} \quad
\ophSigma = O(1),
\end{equation}
which holds, e.g., with high probability when $X$ has independent sub-Gaussian
entries \citep{vershynin2018high}, we can choose $\delta =\log n$ to yield an 
excess risk bound of \smash{$\tilde O(b\sigma\sqrt{d / n})$} with probability at
least $1-1/n$.

Next we give our main result for the ridge-penalized GLM estimator.
This is accomplished by tailoring the general excess risk bound in Proposition 
\ref{prop:risk-ridgeglm} to specific GLMs by identifying the sub-Gaussian
parameter for the noise $\epsilon_i$. Notably, Poisson regression requires a
truncation argument due to the exponential nature of the noise tail, resulting
in a slightly weaker guarantee than that for linear and logistic regression.    

\begin{theorem}
\label{thm:glm-ridge-risk}
Consider the ridge-regularized GLM estimator \smash{$\hthetalambda$} in
\eqref{eq:glm-ridge} with \smash{$\lambda > 0$}. Consider the following cases,
which determine the response distribution for each $i=1,\dots,n$:   
\begin{enumerate}
\item Gaussian (linear regression): \smash{$P_i = \cN(\mu_i,\sigma_i^2)$}.
\item Bernoulli (logistic regression): \smash{$P_i = \mathrm{Bern}(\mu_i)$}.
\item Poisson (Poisson regression): \smash{$P_i = \mathrm{Pois}(\mu_i)$}.
\end{enumerate}
In case 1, we define \smash{$\sigma = \max_{i=1,\dots,n} \sigma_i$}; in case 2,     
\smash{$\sigma = \frac12$}; and in case 3, \smash{$\sigma =
  (6\|\mu\|_\infty+2) \log{n} + 3 \sqrt{\|\mu\|_\infty} $}. Then, for any
\smash{$\delta>0$} and \smash{$b>0$}, choosing $\lambda$ as in
\eqref{eq:glm-ridge-lambda} yields the excess risk bound in
\eqref{eq:glm-ridge-riskhighprob}, which holds with probability at least
\smash{$1 - e^{-\delta}$} in cases 1 and 2, and at least \smash{$1 - e^{-\delta}
  - 1/n$} in case 3 provided $\delta \geq 1$ and $n \geq 3$.  
\end{theorem}

As before, under the spectral assumptions on \smash{$\hSigma$} in
\eqref{eq:hsigma-spectral}, the result in Theorem \ref{thm:glm-ridge-risk}
provides an excess risk bound of \smash{$\tilde O(b\sigma\sqrt{d/n})$} with high 
probability. An important feature of the theorem is its robustness to model
misspecification. To reiterate, the mean $\mu_i$ under the response distribution
$P_i$ need not follow the canonical link of the given GLM. The analysis solely
relies on sub-Gaussian tail bounds for $\epsilon_i$ (and in the Poisson case, an
additional truncation argument). 

As explained above, the intent of Theorem \ref{thm:glm-ridge-risk} is not to
provide the sharpest possible results on excess risk in regularized GLMs. It is
instead to provide a reasonable benchmark to which we can compare risk bounds
for early-stopped gradient descent iterates in the next subsection. Still, it is
worth a brief discussion on how the results in Theorem \ref{thm:glm-ridge-risk} 
compares to what we should expect given the existing literature.

\paragraph{Comparison with existing literature.}

There is of course a huge amount of literature on analyzing risk for regularized
GLMs, and we only provide a brief overview of some of the most relevant points
of comparison. First, taking a broader perspective, an excess prediction risk
rate of \smash{$\sqrt{d/n}$} is standard in statistical learning theory. Such
rates---and generalizations thereof---are to be expected, with the feature
dimension $d$ replaced by a suitable notion of complexity such as VC dimension
\citep{wainwright2019high}.

From the point of view of classical asymptotic theory for maximum likelihood,
one would expect a rate of $d/n$ for estimating GLM parameters in well-specified
settings, though this would be in an asymptotic regime with $d$ fixed and $n \to
\infty$ \citep{wasserman2004all}. Moving to finite-sample prediction risk
analyses for GLMs, in some cases, a rate of $d/n$ is achievable. For example, in
well-specified linear regression without regularization, a standard calculation
yields an expected excess prediction risk of $\sigma^2 d/(2n)$. For logistic 
regression, the story is more nuanced. Using sophisticated analyses that
leverage the self-concordant property of the logistic loss, rates ``close'' to
$d/n$ have been established in a variety of different settings, including
misspecified ones \citep{bach2010self, ostrovskii2021finite}. For Poisson
regression, we are not aware of a finite-sample analysis which produces a rate
of $d/n$ for the excess risk.

\subsection{Risk analysis: early-stopped gradient descent} 
\label{sec:risk-gdglm}

Moving from explicit to implicit regularization, we now study gradient descent
on the unpenalized GLM loss in \eqref{eq:glm-loss}, initialized at $\theta_0=0$.
For simplicity, we focus on a constant step size $\eta_t = \eta$,
$t=0,1,2,\dots$, but the results can be generalized to arbitrary step sizes. We
also consider projected gradient descent on the closed Euclidean ball centered
at the origin \smash{$\ball_d(r) =\{x\in\R^d: \|x\|_2\leq r\}$} of a given
radius $b>0$. 

Toward deriving risk bounds, recall the basic inequality for gradient descent
in Theorem \ref{thm:gd-basicineq}: by \eqref{eq:gd-basicineq}, for the initial
point $\theta_0=0$ and any $\theta\in\R^d$,  
\[
f(\thetagdT) + \frac{\lambda_T}{2} \|\thetagdT - \theta \|_2^2 \leq
f(\theta) + \frac{\lambda_T}{2} \| \theta \|_2^2.
\]
As projected gradient descent is a special case of mirror descent, by Theorem
\ref{thm:md-basicineq} and \eqref{eq:md-basicineq}, the above display is also 
true for this algorithm with $\cC=\ball_d(b)$ and any $\theta\in\ball_d(b)$. 
In either case, this can be rewritten as
\[
f(\thetagdT) - f(\theta) \leq \frac{\lambda_T}{2} \Big( \| \theta \|_2^2 -
\|\thetagdT - \theta \|_2^2 \Big). 
\]
Following the same calculations as for \eqref{eq:glm-risk-diff}, we now add
\smash{$\frac1n \epsilon^\top X (\thetagdT-\theta)$} to each side.
Recalling the decomposition of risk in \eqref{eq:glm-risk}, this gives  
\begin{equation}
\label{eq:glm-risk-diff-gd}
\risk (\thetagdT) - \risk (\theta) \leq \frac1n \epsilon^\top X
\big(\thetagdT -\theta\big) + \frac{\lambda_T}{2} \Big(\|\theta\|_2^2 -
\|\thetagdT - \theta \|_2^2 \Big). 
\end{equation}
The next proposition uses this to derive excess risk bounds which closely
resemble those in Proposition \ref{prop:risk-ridgeglm}. 

\begin{proposition}
\label{prop:risk-gdglm}
Assume that the GLM loss $f$ in \eqref{eq:glm-loss} is $L$-smooth, for
$L>0$. Consider gradient descent \eqref{eq:gd-iterate} with $\eta_t = \eta \in
(0,1/L]$, $t=0,1,2,\dots$, initialized at $\theta_0=0$. Then for any reference
point $\theta \in \R^d$, and for any $T \geq 1$ and $\lambda_T = 1/(\eta T)$,
\[
\risk (\thetagdT) \leq \risk (\theta) + \frac{1}{2\lambda_T}
\Big\|\frac{X^\top\epsilon}{n}\Big\|_2^2 + \frac{\lambda_T}{2}\|\theta\|_2^2. 
\]
If $f$ is only $L$-smooth on $\ball_d(b)$ for some $b>0$, then the same
result holds for projected gradient descent \eqref{eq:pgd-iterate} with $\cC =
\ball_d(b)$, and any reference point $\theta \in \ball_d(b)$.  

Furthermore, assume that $\epsilon_i$, $i=1,\dots,n$, are independent, each
$\epsilon_i \sim \sG(\sigma_i^2)$. Define \smash{$\sigma^2 =
\max_{i=1,\dots,n} \sigma_i^2$} and 
\smash{$\hSigma = \frac1n X^\top X$}. For any $\delta>0$ and $b>0$, define the
target regularization parameter  
\begin{equation} 
\label{eq:glm-gd-lambda}
\lambda^*_T = \frac{\sigma}{b\sqrt{n}} \sqrt{\trhSigma + 2\sqrt\delta\FhSigma +  
  2\delta\ophSigma},
\end{equation}
and define $T^* = 1/(\eta\lambda^*_T)$. If $T^*$ is an integer, then setting 
$T=T^*$ the excess risk of $\thetagdT$ with respect to the class of parameters with
$\ell_2$ norm at most $b$ satisfies 
\begin{equation} 
\label{eq:glm-gd-riskhighprob}
\risk (\thetagdT) - \min_{\theta: \|\theta\|_2\leq b} \risk (\theta) \leq
\frac{b\sigma }{\sqrt{n}} \sqrt{\trhSigma + 2\sqrt\delta\FhSigma +
  2\delta\ophSigma}, 
\end{equation}
with probability at least \smash{$1 - e^{-\delta}$}. If $T^*$ is not an integer,
then setting \smash{$T=\lceil T^* \rceil$}, the above bound holds with an
additional (lower-order) term \smash{$e_T = (\eta \sigma^2 / (2n)) \cdot
  (\trhSigma + 2 \sqrt{\delta} \FhSigma + 2  \delta \ophSigma)$} on the
right-hand side, which is due to discretization.   
\end{proposition}

Under the standard spectral assumptions on \smash{$\hSigma$} in
\eqref{eq:hsigma-spectral}, we can take $\delta =\log n$, and then the bound in   
\eqref{eq:glm-gd-riskhighprob} becomes \smash{$\tilde O(b\sigma\sqrt{d/n})$}, as
before. Moreover, in this case the additional error term (for non-integral
$T^*$) scales as \smash{$e_T = \tilde O(\sigma^2 d/n)$}, confirming that it is
indeed lower-order in the regime $d/n \to 0$.    

Building on this proposition, our next result derives excess risk bounds for
gradient descent for common GLMs. This mirrors Theorem \ref{thm:glm-ridge-risk} 
for ridge-regularized GLMs.

\begin{theorem}
\label{thm:glm-gd-risk}
Let the response distribution $P_i$ and $\sigma$ be defined as in Theorem 
\ref{thm:glm-ridge-risk}, in case 1: Gaussian, case 2: Bernoulli, and case 3: 
Poisson. Note in cases 1 and 2, the GLM loss is $L$-smooth on $\R^d$ with
\smash{$L = \|\hSigma\|_\op$} and \smash{$L = \frac14 \|\hSigma\|_\op$},
respectively; in case 3, it is $L$-smooth on \smash{$\ball_d(b)$} with \smash{$L
  = \exp(b\cdot \max_{i=1,\dots,n} \|x_i\|_2) \cdot \|\hSigma\|_\op$}, for any
$b>0$. Consider gradient descent \eqref{eq:gd-iterate} in cases 1 and 2, and
projected gradient descent \eqref{eq:pgd-iterate} on \smash{$\cC = \ball_d(b)$}
in case 3. In each case, we set $\eta_t = \eta \in (0,1/L]$, $t=0,1,2,\dots$
and initialize at $\theta_0=0$. Fixing any \smash{$\delta>0$} and \smash{$b>0$}, 
let $\lambda_T^*$ be as in \eqref{eq:glm-gd-lambda} and $T^* =
1/(\eta\lambda_T^*)$. The following holds. 
\begin{itemize}
\item If $T^*$ is an integer, then for $T=T^*$ the excess risk bound in
\eqref{eq:glm-gd-riskhighprob} holds with probability at least \smash{$1 
  - e^{-\delta}$} in cases 1 and 2, and at least \smash{$1 - e^{-\delta} - 1/n$}
in case 3 provided that $\delta \geq 1$ and $n \geq 3$.
\item If $T^*$ is not an integer, then for \smash{$T=\lceil T^* \rceil$} the
  excess risk bound holds with the additional (lower-order) additive term $e_T$
  as defined in Proposition \ref{prop:risk-gdglm}.  
\end{itemize}
\end{theorem}

Just like Theorem \ref{thm:glm-ridge-risk}, Theorem \ref{thm:glm-gd-risk}  
translates into a high probability excess risk bound of \smash{$\tilde
  O(b\sigma\sqrt{d/n})$} under the standard spectral assumptions in 
\eqref{eq:hsigma-spectral}. The main point we intend to convey in the above
theorem is that the basic inequality for gradient descent leads to relatively
simple but useful risk bounds for GLMs in general (misspecified) settings,
comparable to those achievable via explicit $\ell_2$ regularization. Below we 
briefly discuss related literature for broader context. 

\paragraph{Comparison with existing literature.} 

For linear regression or quadratic loss more generally there has been a
number of analyses for early-stopped gradient descent/flow that yield
characterizations more precise than the results in the above theorem, such as 
\citet{yao2007early, raskutti2014early, wei2017early, ali2019continuous}. In
particular, the latter authors establish that gradient flow stopped at time $T$
obtains a risk at most 1.69 times that of ridge regression for
$\lambda=1/T$. For logistic regression, time-asymptotic analyses are more
common, which show, e.g., that for linearly separable data gradient descent
converges to the maximum $\ell_2$-margin direction \citep{soudry2018implicit,
  ji2019implicit}. In concurrent work, \citet{wu2025benefits} gave refined risk
bounds for early-stopped gradient descent in logistic regression. This provides
sharper control of the excess risk in certain well-specified settings, at the
rate $d/n$. We are not aware of work analyzing the gradient descent path and its
excess risk for Poisson regression.

\section{Model aggregation with KL regularity}
\label{sec:model-agg}

In this section, we examine mirror descent and its explicit regularization
counterpart, Bregman-divergence-regularization. As our primary application, we 
focus on exponentiated gradient descent and Kullback-Leibler (KL) divergence
regularization. Results for general Bregman divergences are given in the
appendix.  

The notation in this section inherits from the last, with $X \in \R^{n\times d}$ 
denoting a fixed (nonrandom) feature matrix and $Y \in \R^n$ a response vector. 
We again consider the GLM negative log-likelihood function $f$ in 
\eqref{eq:glm-loss}, but now we restrict our attention to parameters $\theta \in
\simplex$, where recall \smash{$\simplex=\{a\in\R^d: a_i\geq 0, \, \sum_{i=1}^d
  a_i = 1\}$} is the probability simplex. As motivation, we may think of this
setup as representing \emph{model aggregation}, with each $X_j$ (the $j\th$
column of $X$) denoting a predictor of the underlying canonical parameter, and
their convex combination \smash{$X\theta = \sum_{j=1}^d \theta_j X_j$}, $\theta
\in \simplex$ representing an aggregate predictor we seek to learn by optimizing
the loss $f$. This approach is often called \emph{linear stacking} in the
literature \citep{wolpert1992stacked, breiman1996stacked}.

The flow of this section is analogous to the last: we first derive simple but
meaningful excess risk bounds for the KL-regularized estimator, which serve as a
benchmark for the early-stopped mirror descent analysis to follow. In each
case, it is not our intention to derive the sharpest results possible, but
instead to demonstrate the power of the basic inequality. 

\subsection{Risk analysis: KL-regularized solution}
\label{sec:risk-klglm}

The KL-regularized GLM estimator is defined by augmenting the GLM loss
\eqref{eq:glm-loss} with a KL penalty, denoted
\begin{equation} 
\label{eq:glm-kl}
\hthetalambda = \argmin_{\theta \in \simplex} \Big( f (\theta) + \lambda \DKL
(\theta, u) \Big),
\end{equation}
where recall \smash{$\DKL(a,b)=\sum_{i=1}^d a_i\log(a_i/b_i)$}, and 
\smash{$u \in \simplex$} is an anchor point. (Note that the
solution exists and is unique because the objective is strictly convex.)  

Using similar arguments to those leading up to \eqref{eq:glm-risk-diff} and
Proposition \ref{prop:risk-ridgeglm} in the ridge regularization setting, the 
proposition below establishes deterministic and high-probability excess risk
bounds.

\begin{proposition}
\label{prop:risk-klglm}
For any \smash{$\lambda>0$}, anchor point \smash{$u\in\simplex$}, and
reference point \smash{$\theta\in\simplex$}, the prediction risk of
\smash{$\hthetalambda$} in \eqref{eq:glm-kl} satisfies   
\[
\risk(\hthetalambda)-\risk(\theta) \leq \frac{1}{\lambda}
\Big\|\frac{X^\top\epsilon}{n}\Big\|_\infty^2 + 2 \lambda \DKL(\theta, u). 
\]
Furthermore, assume that $\epsilon_i$, $i=1,\dots,n$, are independent, each
$\epsilon_i \sim \sG(\sigma_i^2)$, and each \smash{$\|X_j\|_2 \leq \sqrt{n}$}.
Let \smash{$\sigma^2 = \max_{i=1,\dots,n} \sigma_i^2$}. Then, for any
$\delta>0$ and $b>0$, choosing 
\begin{equation}
\label{eq:glm-kl-lambda}
\lambda = \sigma\sqrt{\frac{\log(2d)+\delta}{bn}},
\end{equation}
the excess risk of \smash{$\hthetalambda$} with respect the class of parameters 
with $\DKL(\theta, u) \leq b$ satisfies
\begin{equation}
\label{eq:glm-kl-riskhighprob}
\risk (\hthetalambda) - \min_{\theta: \DKL(\theta, u) \leq b} \risk (\theta)  
\leq 4\sigma\sqrt{\frac{b(\log(2d)+\delta)}{n}},
\end{equation}    
with probability at least \smash{$1- e^{-\delta}$}.
\end{proposition}

Next we specialize this to common GLMs, in a structure similar to Theorem
\ref{thm:glm-ridge-risk} for ridge regularization. 

\begin{theorem}
\label{thm:glm-kl-risk}
Consider the KL-regularized GLM estimator \smash{$\hthetalambda$} in
\eqref{eq:glm-kl} with \smash{$\lambda > 0$}. Consider cases 1, 2, and 3
as described in Theorem \ref{thm:glm-ridge-risk}, and assume that each 
\smash{$\|X_j\|_2 \leq \sqrt{n}$}. Then, for any \smash{$\delta>0$} and
\smash{$b>0$}, choosing $\lambda$ as in \eqref{eq:glm-kl-lambda} yields the
excess risk bound in \eqref{eq:glm-kl-riskhighprob}, which holds with
probability at least \smash{$1 - e^{-\delta}$} in cases 1 and 2, and at least
\smash{$1 - e^{-\delta} - 1/n$} in case 3 provided $\delta \geq 1$ and $n \geq
3$. 
\end{theorem}

Choosing $\delta =\log n$, the bound in \eqref{eq:glm-kl-riskhighprob} becomes
\smash{$\tilde O(\sigma\sqrt{b (\log d)/n})$}, which holds with probability at 
least $1-1/n$ in the Gaussian and Bernoulli cases, and at least $1-2/n$ in the 
Poisson case. We defer discussion of related work until after presenting the
analogous implicit regularization result, in the next subsection. 

\subsection{Risk analysis: early-stopped mirror descent}
\label{sec:risk-egdglm}

We now study early-stopped exponentiated gradient descent as the implicit
regularization counterpart to KL regularization. The proposition below provides
excess risk bounds closely resembling those in Proposition
\ref{prop:risk-klglm}.  

\begin{proposition}
\label{prop:risk-egdglm}
Assume that the GLM loss $f$ in \eqref{eq:glm-loss} is $L$-smooth with
respect to $\|\cdot\|_1$, for $L>0$. Consider exponentiated gradient descent 
\eqref{eq:egd-iterate} with $\eta_t = \eta \in (0,1/L]$, $t=0,1,2,\dots$,
initialized at $\theta_0=u \in \simplex$. Then for any reference
point $\theta \in \simplex$, and for any $T \geq 1$ and $\lambda_T = 1/(\eta T)$, 
\[
\risk (\thetaegdT) - \risk (\theta) \leq \frac{1}{2\lambda_T}
\Big\|\frac{X^\top\epsilon}{n} \Big\|_\infty^2 + \lambda_T \DKL(\theta, u).  
\]
Furthermore, assume that $\epsilon_i$, $i=1,\dots,n$, are independent, each
$\epsilon_i \sim \sG(\sigma_i^2)$, and each \smash{$\|X_j\|_2 \leq \sqrt{n}$}.
Let \smash{$\sigma^2 = \max_{i=1,\dots,n} \sigma_i^2$}. For any $\delta>0$ and
$b>0$, define the target regularization parameter   
\begin{equation} 
\label{eq:glm-egd-lambda}
\lambda_T^*  = \sigma\sqrt{\frac{\log(2d) + \delta}{bn}}.
\end{equation}
Define $T^* = 1/(\eta\lambda_T^*)$. If $T^*$ is an integer, then setting
$T=T^*$ the excess risk of $\thetaegdT$ with respect to the class of parameters
with $\DKL(\theta, u) \leq b$ satisfies 
\begin{equation} 
\label{eq:glm-egd-riskhighprob}
\risk (\thetaegdT) - \min_{\theta: \DKL(\theta, u) \leq b} \risk (\theta) \leq 2
\sigma\sqrt{\frac{b(\log(2d) + \delta)}{n}}, 
\end{equation}
with probability at least \smash{$1 - e^{-\delta}$}. If $T^*$ is not an integer,
then setting \smash{$T=\lceil T^* \rceil$}, the above bound holds with an
additional (lower-order) term \smash{$e_T = \eta \sigma^2 (\log(2d)+\delta) /
  n$} on the right-hand side, which is due to discretization.     
\end{proposition}

Finally, we specialize this to common GLMs, providing results analogous to those
in Theorem \ref{thm:glm-kl-risk}.

\begin{theorem}
\label{thm:glm-egd-risk}
Consider exponentiated gradient descent \eqref{eq:egd-iterate} on the GLM
loss in \eqref{eq:glm-loss}, under cases 1, 2, and 3 as described in Theorem 
\ref{thm:glm-ridge-risk}, and assume that each \smash{$\|X_j\|_2 \leq
  \sqrt{n}$}. These losses are $L$-smooth on $\simplex$ with respect to 
$\|\cdot\|_1$, where the smoothness parameter is \smash{$L = \frac1n
  \sum_{i=1}^n \|x_i\|_\infty^2$} in case 1 (Gaussian), \smash{$L = \frac{1}{4n}
  \sum_{i=1}^n \|x_i\|_\infty^2$} in case 2 (Bernoulli), and \smash{$L = \frac1n
  \sum_{i=1}^n e^{\|x_i\|_\infty} \|x_i\|_\infty^2$} in case 3 (Poisson). 
Set $\eta_t = \eta \in (0,1/L]$, $t=0,1,2,\dots$ and $\theta_0=u \in\simplex$.
Fixing any \smash{$\delta>0$} and \smash{$b>0$}, let
$\lambda_T^*$ be as in \eqref{eq:glm-egd-lambda} and $T^* =
1/(\eta\lambda_T^*)$. The following holds.   
\begin{itemize}
\item If $T^*$ is an integer, then for $T=T^*$ the excess risk bound in
\eqref{eq:glm-egd-riskhighprob} holds with probability at least \smash{$1 
  - e^{-\delta}$} in cases 1 and 2, and at least \smash{$1 - e^{-\delta} - 1/n$}
in case 3 provided that $\delta \geq 1$ and $n \geq 3$.
\item If $T^*$ is not an integer, then for \smash{$T=\lceil T^* \rceil$} the
  excess risk bound holds with the additional (lower-order) additive term $e_T$
  as defined in Proposition \ref{prop:risk-egdglm}.  
\end{itemize}
\end{theorem}

As before, the bound in \eqref{eq:glm-egd-riskhighprob} gives a high probability
excess risk bound of \smash{$\tilde O(\sigma\sqrt{b (\log d)/n})$}. Below we
briefly discuss the broader literature.

\paragraph{Comparison with existing literature.} 

The problem of model aggregation has been studied extensively in the
statistics literature, e.g., \citet{tsybakov2003optimal, juditsky2005recursive,
  lecue2007optimal, juditsky2008learning, dalalyan2009aggregation,
  lecue2013optimality, lecue2014optimal}. A rate of \smash{$\sqrt{(\log d)/n}$}
for excess risk when the comparator class is all convex combinations of base
predictors ($X\theta$ for $\theta \in \simplex$, in our notation) is somewhat
standard. The influential paper \citet{juditsky2005recursive} shows that this
rate is achieved by an estimator called mirror averaging. This performs one epoch
of \emph{stochastic} exponentiated gradient descent followed by an
online-to-batch conversion (in the language of online learning). It is
interesting to emphasize the differences to our results above---we show that 
this rate is achievable using \emph{batch} mirror descent, provided that we use
early stopping.  

It is also worth noting that faster rates are possible under stronger assumptions.
\citet{juditsky2008learning} show that mirror averaging can also produce  
an excess risk bound of order \smash{$(\log d)/n$} when the excess risk is
defined with respect to the single best predictor (rather than the best convex
combination of predictors) in the given finite class. This can also be achieved
by a method called Q-aggregation \citep{lecue2014optimal}, which applies
KL-type regularization to a hybrid loss composed of a linear stacking objective 
(as studied in this section) and a randomized prediction objective (as studied
next).    

\section{Randomized prediction with KL regularity} 
\label{sec:alquier}

We now study \emph{randomized prediction}, an alternative meta-learning approach
to model aggregation as studied previously. While model aggregation outputs a
convex combination of base predictors, randomized prediction samples a single
model according to a learned distribution over the base predictors.  

Using similar notation to the last section, let $X_j \in \R^n$, $j=1,\dots,d$
be base predictors and let $Y \in \R^n$ be a response vector. Given a parameter
$\theta \in \simplex$ on the probability simplex, we define the randomized
predictor $\beta = \beta(X)$ by setting $\beta(X) = X_j$ with probability
$\theta_j$. As shorthand, we will write the corresponding predictor as $\beta
\sim \theta$. We define the training risk as
\[
\hR(\beta) = \frac1n \sum_{i=1}^n \ell(\beta(X)_i, Y_i)
\]
where $\ell$ is an arbitrary loss function. 

A randomized predictor can also make test predictions, and hence we can also 
define a suitable notion of test risk. Let $X_j^* \in \R^n$, $j=1,\dots,d$ and
$Y^* \in \R^n$ denote test instances of the base predictors and response
variable, respectively. Then the randomized predictor $\beta \sim \theta$ sets
$\beta(X^*) = X_j^*$ with probability $\theta_j$. We define the test risk as
\[
R(\beta) = \E\bigg[ \frac1n \sum_{i=1}^n \ell(\beta(X^*)_i, Y^*_i) \bigg],  
\]
where the expectation is with respect to any joint distribution $P$ on
$(X^*,Y^*)$. This setting is quite general. For $\ell$ taking the GLM form, and 
the joint distribution $P$ placing a point mass at $X^*=X$, this recovers the
setup in the previous two sections, with $R(\beta)=\risk(\theta_j)$ and
\smash{$\hR(\beta)=f(\theta_j)$} for a particular (random) index $j$. However,
the current setup is even broader, allowing $X^*$ to be random, and $\ell$ to be
arbitrary.  
 
To construct the sampling distribution, we will consider exponential weighting
based on the empirical risk evaluated at each base model:  
\begin{equation} 
\label{eq:gibbs-formula}
\hthetalambda(\de \beta) \propto \exp(- \hR (\beta) / \lambda) \cdot 
u (\de \beta), 
\end{equation}
where $u\in\simplex$ is a prior measure and
\smash{$\lambda>0$} is a tuning parameter. With no prior information, one may
take the uniform \smash{$u = \pi$}
on $\simplex$. The above estimator is called the Gibbs posterior in the Bayesian
statistics literature. It can be equivalently defined by solving a
KL-regularized optimization problem \citep{alquier2024user}:    
\begin{equation}
\label{eq:gibbs-min}
\hthetalambda = \argmin_{\theta \in \simplex} \Big( \E_{\beta \sim \theta}
[\hR(\beta)] + \lambda \DKL (\theta, u) \Big), 
\end{equation}
where the expectation is only with respect to the randomness in drawing $\beta
\sim \theta$. This problem is also referred to as information risk minimization
in some parts of the literature \citep{leung2006information,
  zhang2006information, xu2017information}.   

The flow of the previous two sections would suggest that we would here produce 
a risk analysis of the explicitly-regularized estimator \eqref{eq:gibbs-min} and 
then compare to what is possible to derive from the basic inequality for its 
implicitly-regularized counterpart, exponentiated gradient descent on the loss 
\smash{$f(\theta) = \E_{\beta \sim \theta} [\hR(\beta)]$}. In this section,
however, we depart from this strategy for one important reason: these two 
estimators actually coincide exactly. 

\begin{proposition}
\label{prop:gibbs-equiv}
Consider the Gibbs posterior \smash{$\hthetalambda$} in \eqref{eq:gibbs-min}
with a prior $u\in\simplex$, and exponentiated gradient descent 
\eqref{eq:egd-iterate} with linear loss \smash{$f(\theta) = \E_{\beta \sim
    \theta} [\hR(\beta)]$}, constant step sizes $\eta_t = \eta > 0$,
$t=0,1,2,\dots$, and initialization $\theta_0=u$. Then for \smash{$\lambda =
  1/(\eta T)$}, these two estimators coincide: \smash{$\hthetalambda = \thetaegdT$}.   
\end{proposition}

This means that to derive excess risk bounds, we can apply basic inequality 
arguments either to \smash{$\hthetalambda$} or to $\thetaegdT$. In what follows,
we pursue the latter, because it results in sharper constants, i.e., the excess
risk bounds are sharper by a factor of 2 (which is unsurprising, since this was
also the case in the last two sections). We restrict our attention to the case
in which the optimal tuning corresponds to an integral time $T$, for brevity
(and analogous discretization results would follow as before).

\begin{proposition}
\label{prop:gibbs-risk}
For any \smash{$\lambda>0$}, prior \smash{$u\in\simplex$}, and reference point
\smash{$\theta\in\simplex$}, the risk of \smash{$\hthetalambda$} in
\eqref{eq:gibbs-min} satisfies    
\[
\E_{\beta \sim \htheta_\lambda} [R(\beta)] - \E_{\beta \sim \theta} [R(\beta)]
\leq \frac{1}{2\lambda} \|\hR - R\|_\infty^2 + \lambda \DKL (\theta, u),
\]
where \smash{$\|\hR - R\|_\infty = \sup_\beta |\hR(\beta)-R(\beta)|$}. The same
bound applies to exponentiated gradient descent \eqref{eq:egd-iterate} with the
loss \smash{$f(\theta) = \E_{\beta \sim   \theta} [\hR(\beta)]$} (and constant
step size $\eta$ and $\theta_0=u$), provided $T=1/(\eta\lambda)$ is an integer.  

Furthermore, if the loss $\ell$ is bounded by \smash{$C>0$} and the training
samples are i.i.d., then for any \smash{$\delta>0$} and \smash{$b>0$}, it holds
that 
\[
\E_{\beta \sim \htheta_\lambda} [R(\beta)] - 
\min_{\theta: \DKL(\theta, u) \leq b} \E_{\beta \sim \theta} [R(\beta)] 
\leq \frac{C^2 (\log(2d) + \delta)}{4n \lambda} +  \lambda b.
\]
with probability at least \smash{$1-e^{-\delta}$}.	In particular, choosing 
\[
\lambda = \frac{C}{2}\sqrt{\frac{\log(2d)+\delta}{bn}},
\]
it holds that
\[
\E_{\beta \sim \htheta_\lambda} [R(\beta)] - 
\min_{\theta: \DKL(\theta, u) \leq b} \E_{\beta \sim \theta} [R(\beta)] 
\leq C \sqrt{\frac{b(\log(2d)+\delta)}{n}}.
\]
with probability at least \smash{$1-e^{-\delta}$}. The same result again applies
to exponentiated gradient descent provided $T=1/(\eta\lambda)$
is an integer. 
\end{proposition}

Lastly, we briefly discuss related literature.

\paragraph{Comparison with existing literature.}

The most directly related analysis given by \citet{alquier2024user}; under the
same conditions as our result, this author establishes that the Gibbs posterior
estimator satisfies   
\[
\E_{\beta\sim\hthetalambda}[R(\beta)] -
\min_{\theta: \DKL(\theta, u) \leq b} \E_{\beta \sim \theta} [R(\beta)]  
\leq \frac{C^2}{4n\lambda} + 2\lambda ( b + \log 2 + \delta),
\]
with probability at least \smash{$1-e^{-\delta}$}. Compared to the first result
in Proposition \ref{prop:gibbs-risk}, we can see that this is sharper by a
factor of $\log d$ in the first term on the right-hand side (and so in
particular, it applies in situations where the number of base models is
infinite). This is due to a more sophisticated analysis of the empirical process
term in \citet{alquier2024user}; they use the Donsker-Varadhan formula, whereas
we simply use Young's inequality and Hoeffding's inequality. Still, it is
interesting to see that the basic inequality for mirror descent produces
an excess risk rate of \smash{$\sqrt{(\log d)/n}$}, which is a meaningful result
for randomized prediction and only slightly worse, for finite $d$, than
the rate \smash{$1/\sqrt{n}$} achievable via more advanced techniques.

\section{Experiments}
\label{sec:experiments}

We present empirical results to supplement our theoretical findings on the 
relationship between explicit and implicit regularization, investigating both
training dynamics and prediction risk. Python code to reproduce our experiments
is available at \url{https://github.com/100shpaik/}.   

\paragraph{Experimental setup.}

We evaluate gradient descent (GD) and exponentiated gradient descent (EGD), as
in \eqref{eq:gd-iterate} and \eqref{eq:egd-iterate}, respectively. In all
experiments that follow, we initialize GD at the origin, and EGD at the
uniform distribution. We apply GD to three GLM losses: linear, logistic, and
Poisson (as in Section \ref{sec:glm}), and compare its iterates to explicit
ridge regularization \eqref{eq:glm-ridge}. We run EGD on these same losses,
in a model aggregation context (as in Section \ref{sec:model-agg}), and compare
its iterates to explicit KL regularization \eqref{eq:glm-kl}. In general, we use
variable step sizes \smash{$\eta_t$}, $t=0,1,2,\dots$, and define the total
elapsed time as \smash{$\tau = \sum_{t=0}^{T-1} \eta_t$}, where we align
$1/\tau$ with the regularization parameter $\lambda$. Further implementation
details are provided in the appendix.

\paragraph{Data distributions.}

We study both underparameterized (\smash{$n>d$}) and overparameterized
(\smash{$n<d$}) regimes, and generate features \smash{$X \in \R^{n \times d}$}
and responses \smash{$Y \in \R^n$} as follows. We sample the entries of $X$
independently from \smash{$\cN(0,1)$}, and then sample $Y$ given $X$ from a 
well-specified model for each GLM, with parameter $\theta\in\R^d$. In the
regression tasks (where we compare GD and ridge regularization), we sample
the entries of $\theta$ uniformly on $[-1,1]$. In model aggregation tasks (where
we compare EGD and KL regularization), we sample uniformly on $[0,1]$, and then
renormalize so that $\theta\in\simplex$.    

We introduce an additional parameter \smash{$\gamma>0$} to control the
signal-to-noise ratio. To be precise, for each $i=1,\dots,n$, we independently
sample $Y_i|x_i$ in the linear, logistic, and Poisson regression cases as
follows:      
\begin{itemize}
\item linear: \smash{$Y_i|x_i \sim \cN(x_i^\top \theta, \gamma^2)$}; 
\item logistic: \smash{$Y_i|x_i \sim \mathrm{Bernoulli}(p_i)$} with
    \smash{$p_i=1/(1+\exp(-\gamma x_i^\top \theta))$};
\item Poisson: \smash{$Y_i|x_i \sim \mathrm{Pois}(\mu_i)$} with 
  \smash{$\mu_i = \exp(\gamma x_i^\top \theta)$}. 
\end{itemize}
Table \ref{table:experiment-setting} summarizes the values of $n$, $d$, and
$\gamma$ used in our experiments. (The specific values of $\gamma$ are chosen to
produce ``U-shaped'' prediction risk curves, avoiding regimes where the risk is
uninterestingly monotonic, which typically occurs when $\gamma$ is too small or
large.) 

\begin{table}[htb]
\centering
\caption{Summary of $n$, $d$, and $\gamma$ values used in our experiments.}
\begin{tabular}{| l | c | c | c | c |}
\hline
& \multicolumn{2}{c|}{\textbf{Regression (GD vs.\ ridge)}} &
\multicolumn{2}{c|}{\textbf{Model aggregation (EGD vs.\ KL)}} \\ 
\cline{2-5}
\textbf{GLM} & Underparameterized & Overparameterized & Underparameterized &
  Overparameterized \\ 
& $(n,d) = (200,20)$ & $(n,d) = (100, 200)$ & $(n,d) = (200,20)$ & $(n,d) = (30, 60)$ \\
\hline
Linear & $\gamma = 5.0$ & $\gamma = 5.0$ & $\gamma = 1.0$ & $\gamma = 0.1$ \\
Logistic & $\gamma = 0.3$ & $\gamma = 0.5$ & $\gamma = 1.5$ & $\gamma = 10.0$ \\
Poisson & $\gamma = 0.1$ & $\gamma = 0.15$ & $\gamma = 1.2$ & $\gamma = 3.5$ \\
\hline
\end{tabular}
\label{table:experiment-setting}
\end{table}

\subsection{Training dynamics}

\begin{figure}[htb]
\centering
\begin{subfigure}[b]{0.48\textwidth}
\centering
\includegraphics[width=\textwidth]{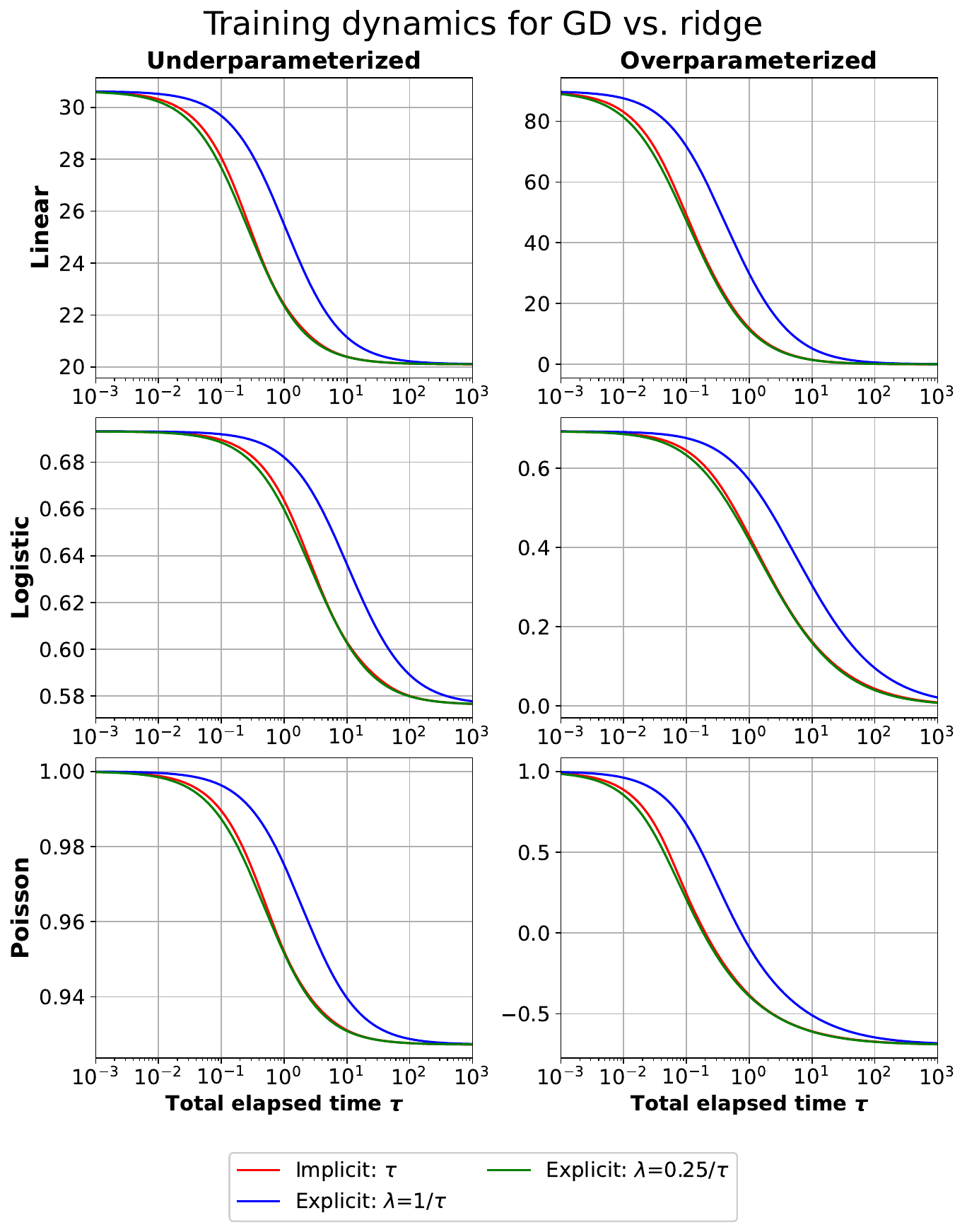}
\caption{GD vs.\ ridge, plotted with $\|\theta\|_2^2$ penalty.}
\end{subfigure}
\hfill
\begin{subfigure}[b]{0.48\textwidth}
\centering
\includegraphics[width=\textwidth]{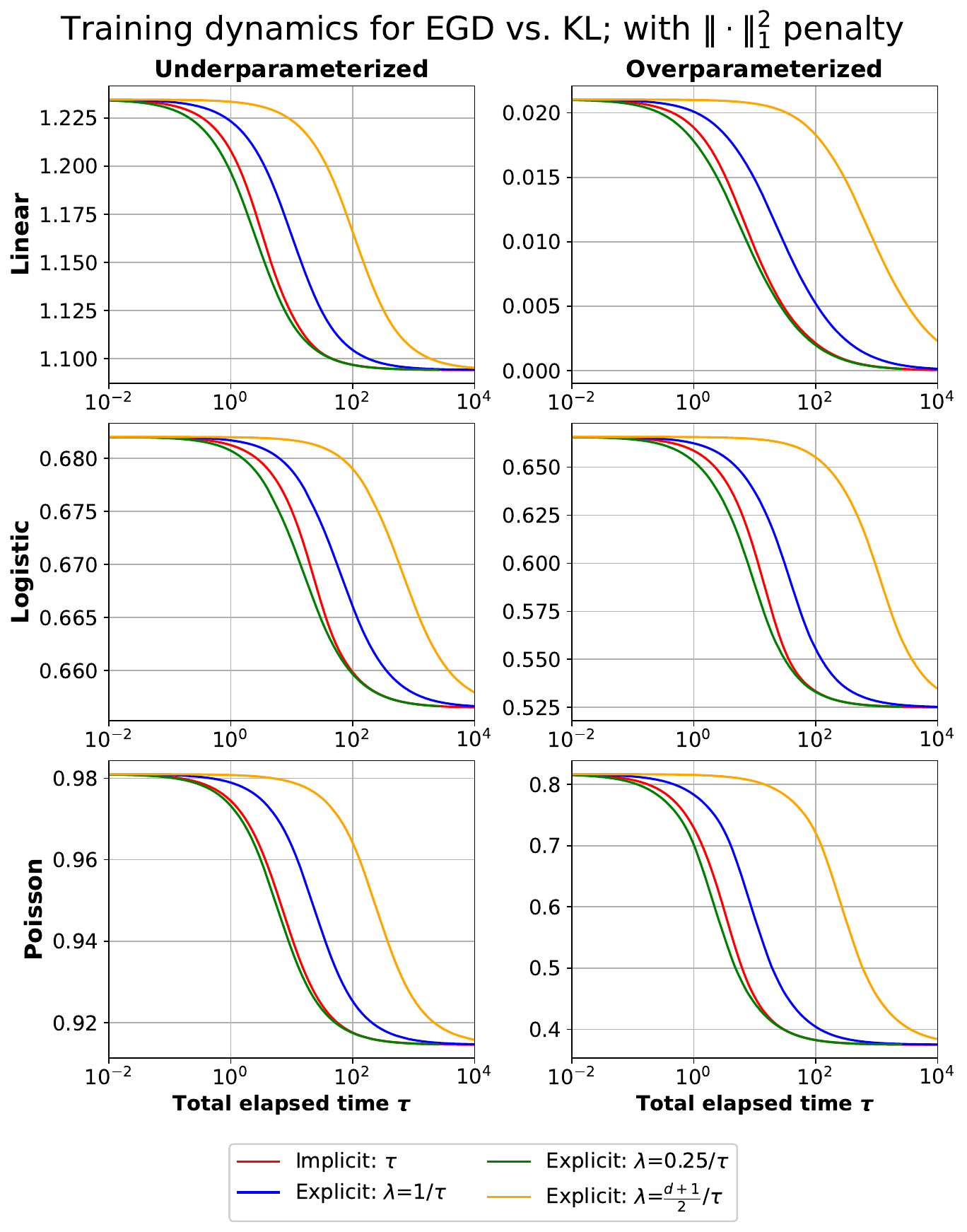}
\caption{EGD vs.\ KL, plotted with $\|\theta-\pi\|_1^2$ penalty.}
\end{subfigure}
\caption{Training envelopes for (a) regression and (b) model aggregation tasks.} 
\label{fig:training-loss-penalty-trajectory}
\end{figure}

Figure \ref{fig:training-loss-penalty-trajectory} plots the penalized
training loss: 
\[
P(\theta,\lambda) = f(\theta) + \lambda g(\theta), 
\]
across implicitly- and explicitly-regularized estimates. In panel (a), we
consider the regression tasks, where $g(\theta) = \|\theta\|_2^2$, and in panel
(b), we consider model aggregation, where $g(\theta) = \|\theta-\pi\|_1^2$. Each
row corresponds to a different GLM loss $f$: linear, logistic, and Poisson. 

Examining first the results in panel (a), we compare three curves: one for GD  
iterates, \smash{$P(\thetagdT, \frac{1}{4\tau})$} as the total elapsed time
$\tau$ varies (plotted in red), and two for ridge estimates,
\smash{$P(\hthetalambda, \lambda)$} as $\lambda$ varies, using \smash{$\lambda 
  = \frac{1}{4\tau}$} (in green) or  \smash{$\lambda = \frac{1}{\tau}$} (in
blue). According to Corollary \ref{cor:gd-envelope}, we know that the red curve
must lie in between the green and blue, and we see this is indeed true in the
figure. Notably, the red line closely tracks the green line for
\smash{$\lambda=\frac{1}{4\tau}$}. This reveals a stronger correspondence
between $\thetagdT$ and \smash{$\hthetalambda$} with
\smash{$\lambda=\frac{1}{4\tau}$} than anticipated by the existing theory,
which will be revisited in the prediction risk analysis shortly.  

The results in panel (b) are analogous, except with squared $\ell_1$ penalties
used in $P(\theta,\lambda)$, as suggested by Corollary
\ref{cor:md-envelope}. Here we have an additional curve corresponding to
explicit (KL) regularization, where \smash{$\lambda = \frac{d+1}{2\tau}$} (as
in the upper bound in the corollary, recalling $G=d$ and $\alpha=1$ in the
present setting). We can see that this is overly conservative in practice, and
the correspondence is tightest between $\thetaegdT$ and \smash{$\hthetalambda$}
with \smash{$\lambda=\frac{1}{4\tau}$}. Moreover, results where
$P(\theta,\lambda)$ is defined in terms of a KL penalty (not covered by any
existing theory) are similar, and given in the appendix. 

\subsection{Prediction risk}

\begin{figure}[htb]
\centering
\begin{subfigure}[b]{0.48\textwidth}
\centering
\includegraphics[width=\textwidth]{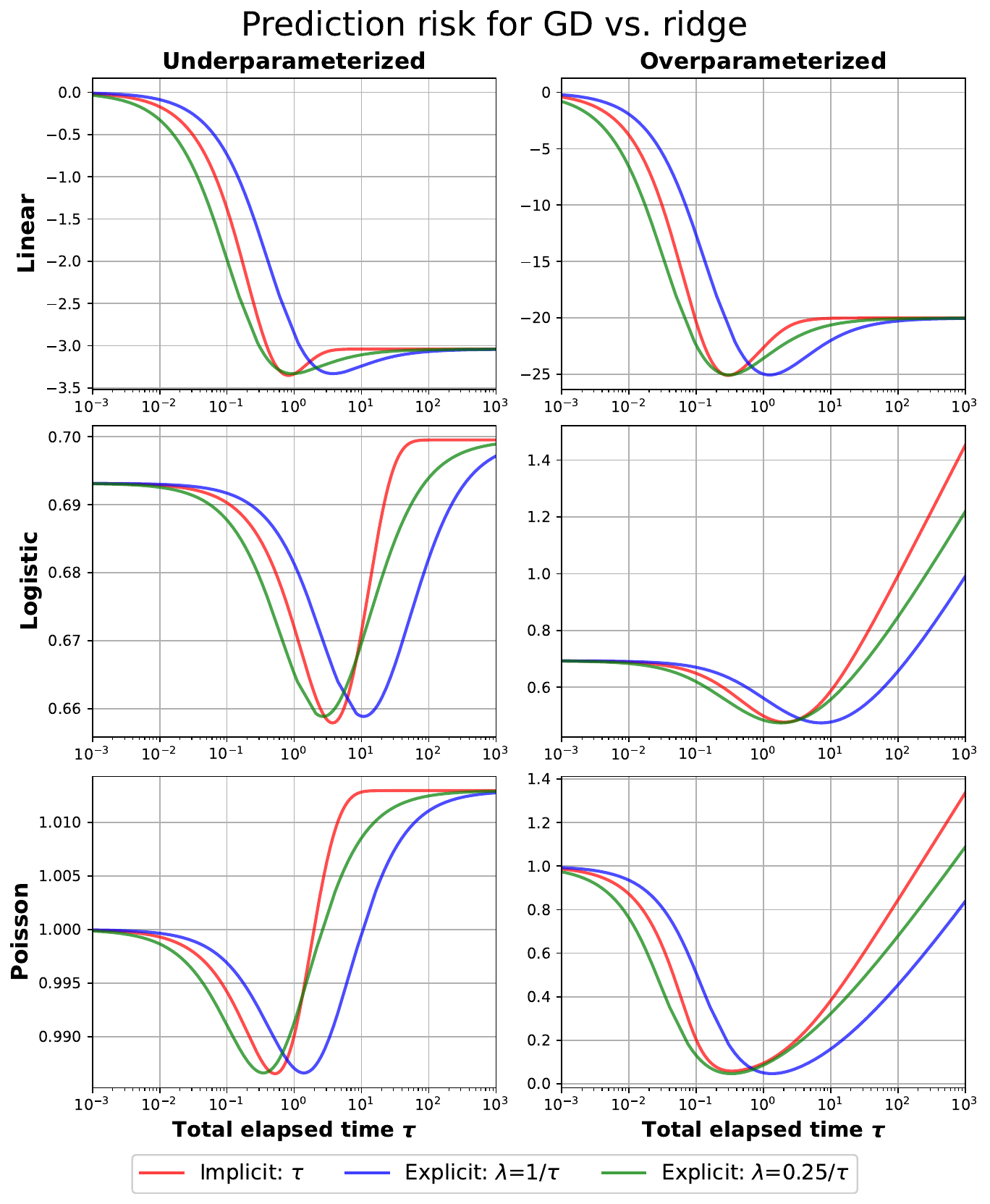}
\caption{GD vs.\ ridge}
\end{subfigure}
\hfill
\begin{subfigure}[b]{0.48\textwidth}
\centering
\includegraphics[width=\textwidth]{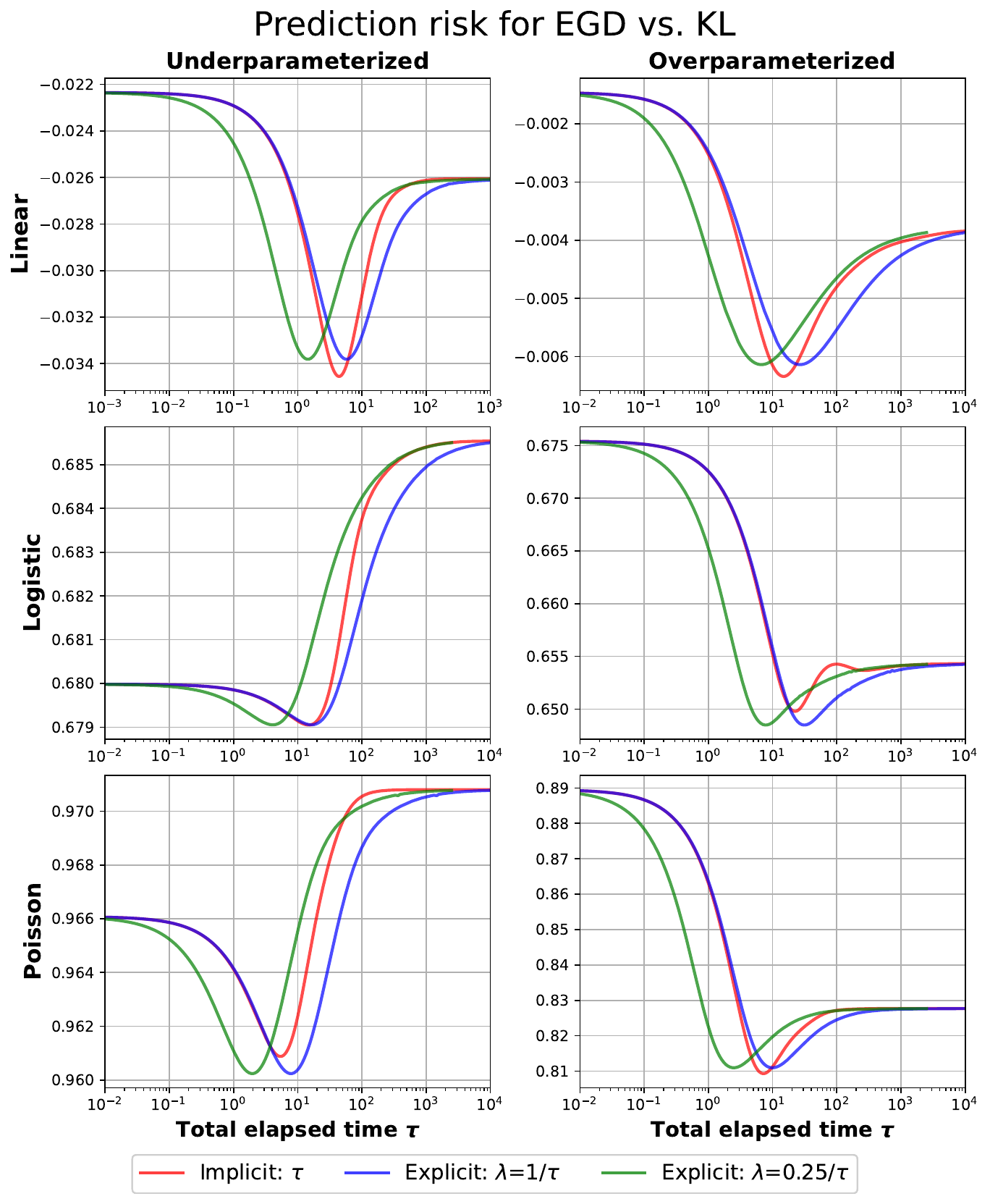}
\caption{EGD vs.\ KL}
\end{subfigure}
\caption{Prediction risk curves for (a) regression and (b) model aggregation
  tasks.}  
\label{fig:risk-comparison}
\end{figure}

Figure \ref{fig:risk-comparison} plots prediction risk curves (where
this is defined in a fixed-X sense, as explained in Section
\ref{sec:glm}) in a format analogous to that used in Figure
\ref{fig:training-loss-penalty-trajectory}. Recall that unlike our results on 
training dynamics, the risk theory in our paper does not draw direct
comparisons between the curves for implicitly- and explicitly regularized
estimates; instead, Theorems \ref{thm:glm-ridge-risk}--\ref{thm:glm-egd-risk}
derive excess risk bounds, which end up being very similar for GD and ridge
regularization, and EGD and KL regularization. Nevertheless, we can see in panel
(a) that the risk curves for GD and ridge estimates are qualitatively quite
similar, with the red and green curves---which correspond to $\thetagdT$ and 
\smash{$\hthetalambda$} with \smash{$\lambda=\frac{1}{4\tau}$}---being overall
the closest across the settings. This is also broadly true in panel (b), for the
EGD and KL-regularized estimates, but now the red and blue
curves---which correspond to $\thetaegdT$ and \smash{$\hthetalambda$} with  
\smash{$\lambda=\frac{1}{\tau}$}---being closer. Finally, the minimum risk
obtained by implicit and explicit regularization is quite close in all cases, 
with neither one dominating the other.

\subsection{Solution path}

\begin{figure}[p]
\centering
\begin{subfigure}[b]{0.8\textwidth}
\centering
\includegraphics[width=\textwidth]{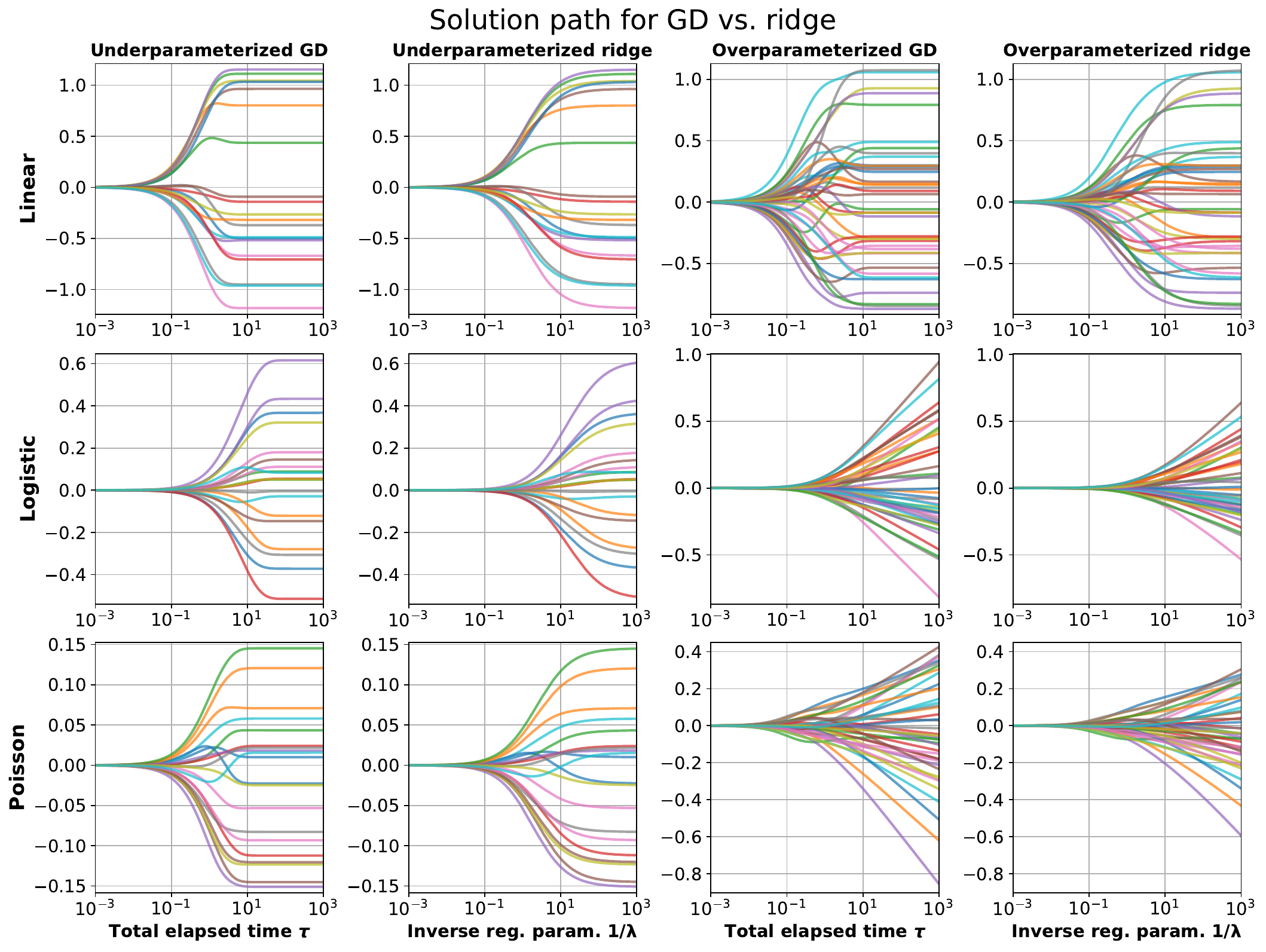}
\caption{GD vs.\ ridge}
\end{subfigure}

\bigskip
\begin{subfigure}[b]{0.8\textwidth}
\centering
\includegraphics[width=\textwidth]{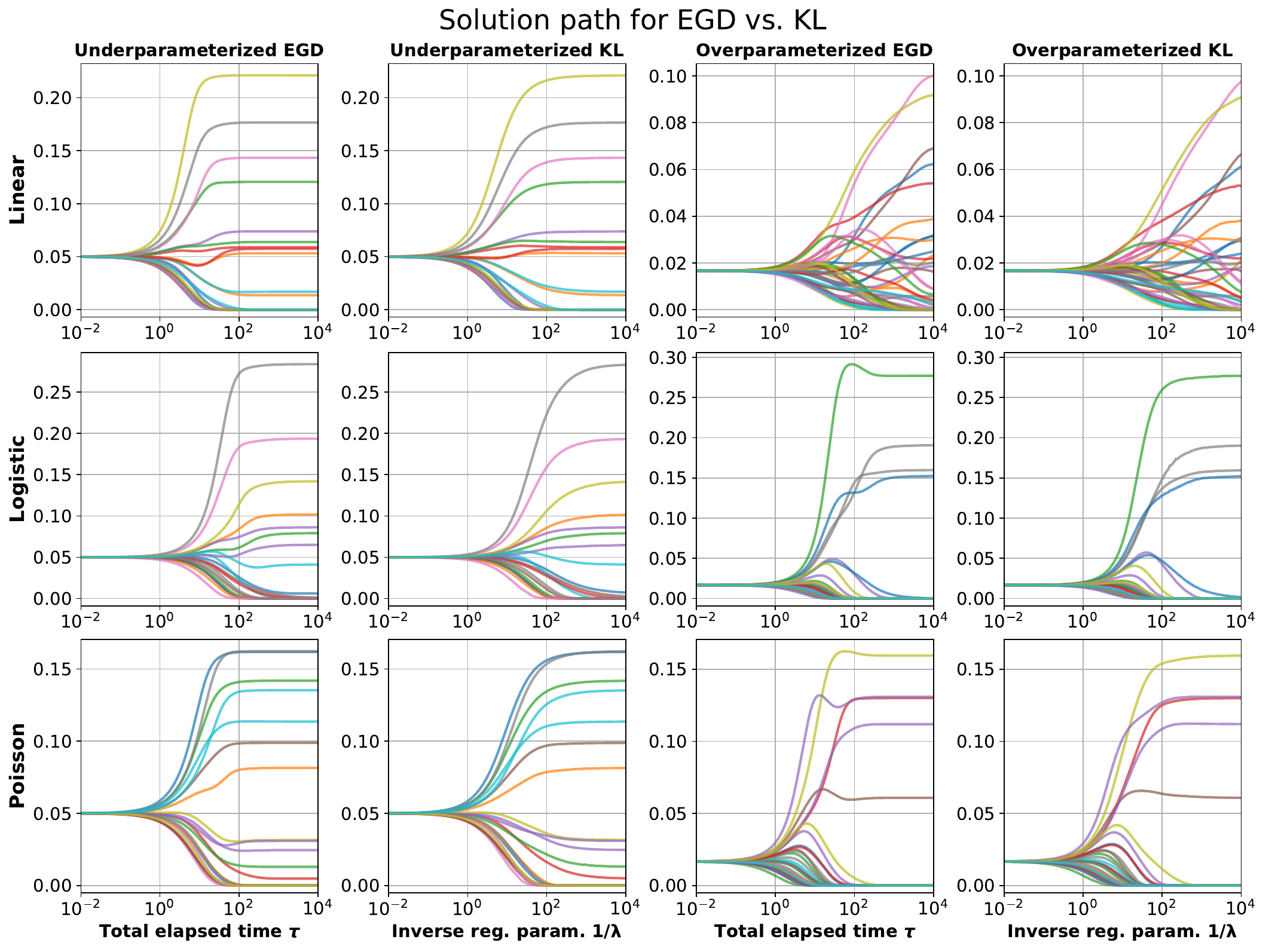}
\caption{EGD vs.\ KL}
\end{subfigure}
\caption{Solution paths for (a) regression and (b) model aggregation. We plot
  all \smash{$d=20$} coordinate paths for underparameterized regimes and the
  first $40$ coordinates for overparameterized regimes, to avoid overcrowding
  the visualization.}
\label{fig:training-sol-path}
\end{figure}

Figure \ref{fig:training-sol-path} displays solution paths for implicitly- and
explicitly-regularized estimates. In a given plot, one curve represents one
coordinate path, either $\theta_{Ti}$ for varying $\tau$ or
\smash{$\htheta_{\lambda i}$} for varying $\lambda$. Although none of our theory
in this paper compares solution paths directly, we generally see a striking
resemblance between GD paths and ridge-regularized paths, and between EGD paths
and KL-regularized paths. This provides further compelling evidence of the deep
connections between implicit and explicit regularization. Meanwhile, it is
interesting to note that the paths for EGD and KL regularization converge to a
sparse solution in the limit, as $\tau \to \infty$ or $\lambda \to 0$. Studying
connections between these regularization mechanisms and lasso ($\ell_1$)
regularization may be a topic of future work.

\clearpage
\section{Other basic inequalities}
\label{sec:other-algos}

We also discuss basic inequalities for two additional first-order iterative
algorithms, proximal gradient descent and NoLips (essentially, a version of
mirror descent which uses different assumptions). 

\subsection{Proximal gradient descent}

Proximal gradient descent applies to composite objectives \smash{$f=f_0+f_1$}
with $f_0,f_1$ convex and $f_0$ differentiable. Given an initial point
$\theta_0\in\R^d$, and step sizes \smash{$\{\eta_t\}_{t=0}^\infty$}, this method
generates iterates for $t=0,1,2,\dots$ via  
\begin{equation} 
\label{eq:proxgd-iterate}
\theta_{t+1} = \prox_{\eta_t f_1}\Big( \theta_t - \eta_t \nabla f_0(\theta_t)
\Big),
\end{equation}
where the proximal operator is defined as
\[
\prox_{\eta f_1}(\theta) = \argmin_{z\in\R^d} \bigg( \frac12 \|\theta-z\|_2^2 +  
\eta f_1(z) \bigg).
\]
A basic requirement for using this algorithm is that this proximal operator can
be computed in closed-form or efficiently approximated. The following theorem
derives a basic inequality for proximal gradient descent, which generalizes
Theorem \ref{thm:gd-basicineq} (which corresponds to $f_1=0$). 

\begin{theorem}
\label{thm:proxgd-basicineq}
Let $f_0:\R^d\to \R$ be convex, differentiable, and $L$-smooth for
$L>0$, and $f_1:\R^d\to \R$ be convex. Consider proximal gradient descent
\eqref{eq:proxgd-iterate} with step sizes $\eta_t\in (0,1/L]$,
$t=0,1,2,\dots$. For any reference point $z\in\R^d$, and any iteration $T \geq
1$, it holds for $f=f_0+f_1$ that
\[
f(\theta_T) - f(z) \leq \frac{1}{2\sum_{t=0}^{T-1} \eta_t} \Big( \|\theta_0 -
z\|_2^2 - \|\theta_T - z\|_2^2 \Big).
\]
\end{theorem}

By the same arguments used to derive training envelopes in Corollary
\ref{cor:gd-envelope}, the result in Theorem \ref{thm:proxgd-basicineq} implies    
\begin{multline*}
\min_{z\in\R^d} \Big( f_0(z) + f_1(z) + \frac{\lambda_T}{4} \|\theta_0 - z\|_2^2
\Big) \leq f_0(\theta_T) + f_1(\theta_T) + \frac{\lambda_T}{4} \|\theta_0 -
\theta_T\|_2^2 \\ \leq \min_{z\in\R^d} \Big( f_0(z) + f_1(z) + \lambda_T
\|\theta_0 - z\|_2^2 \Big).   
\end{multline*}
When $f_1(\theta) = \lambda \|\theta\|_1$ is a lasso penalty, the corresponding
proximal operator is known as soft-thresholding, and proximal gradient descent
is called the iterative soft-thresholding algorithm (ISTA). In this case,
further specifying, e.g., \smash{$f_0(\theta) = \frac{1}{2n} \|Y-X\theta\|_2^2$}
and initializing $\theta_0=0$, the above display becomes      
\begin{multline*}
\min_{z\in\R^d} \bigg( \frac{1}{2n} \|Y-Xz\|_2^2 + \lambda \|z\|_1 +
\frac{\lambda_T}{4} \|z\|_2^2 \bigg) \leq \frac{1}{2n} \|Y-X\theta_T\|_2^2 + 
\lambda \|\theta_T\|_1 + \frac{\lambda_T}{4} \|\theta_T\|_2^2 \\ 
\leq \min_{z\in\R^d} \bigg( \frac{1}{2n} \|Y-Xz\|_2^2 + \lambda \|z\|_1 + 
\lambda_T \|z\|_2^2 \bigg).
\end{multline*}
The interpretation is that early-stopped ISTA provides regularization akin to
explicit elastic net (composite lasso and ridge penalties) regularization. 

\subsection{NoLips algorithm}

\citet{bauschke2017descent} proposed an iterative algorithm called NoLips, which
can be viewed as an instance of mirror descent but operates under nonstandard 
assumptions; notably, it does not require Lipschitz smoothness of the objective
$f$ or strong convexity of the potential function $\phi$. Below we derive a
basic inequality for this algorithm.  

\begin{theorem}
\label{thm:nolips-basicineq}
Let $\Omega,\cC$ be closed, convex sets in $\R^d$ with nonempty interior, and 
$\cC\subseteq\Omega$. Assume $f:\Omega\to\R$ is convex on $\cC$, and
differentiable on \smash{$\interior(\Omega)$}. Assume
$\phi:\Omega\to\R$ is of Legendre type. Moreover, assume that $L\phi-f$ is
convex on \smash{$\cC\cap\interior(\Omega)$}, for some $L>0$. Consider
mirror descent \eqref{eq:md-iterate}---which in this case, is known as
NoLips---initialized at $\theta_0\in\cC\cap\interior(\Omega)$ with step
sizes $\eta_t\in(0,1/L]$, $t=0,1,2,\dots$. For any reference point $z\in\cC$,
and any $T \geq 1$,       
\[
f(\theta_T) - f(z) \leq \frac{1}{\sum_{t=0}^{T-1}\eta_t}\Big(\Dphi(z, \theta_0)
- \Dphi(z, \theta_T)\Big). 
\]
\end{theorem}

An interesting example where the assumptions of Theorem
\ref{thm:nolips-basicineq} hold but those of Theorem \ref{thm:md-basicineq} (the 
standard mirror descent assumptions) do not hold is the Poisson linear inverse (PLI)  
problem. Unlike a Poisson GLM which uses the canonical link \smash{$\log(\mu) =
  X\theta$}, in PLI we use a linear link \smash{$\mu = X\theta$}. This is
typically done in situations where the elements of $X$ are nonnegative and we
constrain those of $\theta$ to be nonnegative, ensuring that $\mu$ itself has
nonnegative entries. The Poisson negative log-likelihood under this link
function is (after scaling by \smash{$\frac1n$}):  
\[
f(\theta) = \frac1n \sum_{i=1}^n \Big( -Y_i\log(x_i^\top \theta) + x_i^\top
\theta \Big).  
\]
This is not Lipschitz smooth over the nonnegative orthant $\R_+^d$, because
the derivative of $\log(u)$ blows up as $u \to 0^+$. This means that the
traditional mirror descent analysis cannot be applied. However, as shown in 
\citet{bauschke2017descent}, if we choose Burg's entropy as the potential
function \smash{$\phi(\theta) = -\sum_{i=1}^d \log(\theta_i)$}, then one can
show that $L\phi-f$ is convex on $\R_{++}^d$ for any $L \geq \|Y\|_1$. Thus
Theorem \ref{thm:nolips-basicineq} applies, and the interpretation loosely
speaking is that early-stopping the NoLips algorithm for the PLI problem
provides regularization akin to the Bregman divergence of Burg's entropy: this 
is \smash{$\Dphi(\theta, \theta_0) = \sum_{i=1}^d \frac{\theta_i}{\theta_{0i}} -
  \log(\frac{\theta_i}{\theta_{0i}})$} (up to constants), and is known as the
Itakura-Saito divergence.  

\section{Discussion}
\label{sec:discussion} 

In this paper, we examined basic inequalities for first-order optimization
algorithms, which provide a unified framework for analyzing implicit
regularization through both computational and statistical perspectives. We
demonstrated the broad applicability of this framework by showing it can be used
to infer properties about training dynamics and also to derive excess risk
bounds, for different algorithms in various settings.

We close by highlighting a few directions for future work. The first is to
derive refined basic inequalities under stronger assumptions, such as strong
convexity. We believe this may enable us to derive faster excess risk rates for
early-stopped gradient or mirror descent (i.e., of order $d/n$ for GLM
regression, or $\log{d}/n$ for model aggregation). A second direction of
interest is to extend the basic inequality framework to stochastic gradient
methods, ideally permitting some degree of nonconvexity in the objective. A
third direction would be to establish basic inequalities for steepest descent
algorithms with respect to non-Euclidean norms, such as the $\ell_1$ norm, in
which case steepest descent becomes greedy coordinate descent, also called
forward stagewise regression. Given the rich tradition of work connecting
stagewise algorithms to $\ell_1$ regularization in statistics
\citep{efron2004least}, it may be interesting to revisit this problem and see if
a perspective based on something like the basic inequality can shed new light.

\subsection*{Acknowledgements}

SP and RJT were supported by Office of Naval Research (ONR) grant number
N00014-20-1-2787. KZ was supported by the Founder’s Postdoctoral Fellowship in  
Statistics at Columbia University.

{\RaggedRight
\bibliographystyle{plainnat}
\bibliography{main.bib}}

@book{buhlmann2011statistics,
  author = {Peter Buhlmann and Sara van de Geer},
  publisher = {Springer},
  title = {Statistics for High-Dimensional Data},
  year = {2011}}

@book{nesterov2003introductory,
  title={Introductory Lectures on Convex Optimization: A Basic Course},
  author={Nesterov, Yurii},
  year={2003},
  publisher={Springer}
}

@article{nemirovski2009robust,
  title={Robust stochastic approximation approach to stochastic programming},
  author={Nemirovski, Arkadi and Juditsky, Anatoli and Lan, Guanghui and Shapiro, Alexander},
  journal={SIAM Journal on Optimization},
  volume={19},
  number={4},
  pages={1574--1609},
  year={2009}
}

@inproceedings{reddi2018convergence,
  title={On the Convergence of {Adam} and Beyond},
  author={Reddi, Sashank J. and Kale, Satyen and Kumar, Sanjiv},
  booktitle={International Conference on Learning Representations},
  year={2018}
}

@article{wu2025benefits,
  title={Benefits of Early Stopping in Gradient Descent for Overparameterized Logistic Regression},
  author={Wu, Jingfeng and Bartlett, Peter L. and Telgarsky, Matus and Yu, Bin},
  journal={arXiv preprint arXiv:2502.13283},
  year={2025}
}

@inproceedings{wu2024large,
  title={Large stepsize gradient descent for logistic loss: {Non-monotonicity} of the loss improves optimization efficiency},
  author={Wu, Jingfeng and Bartlett, Peter L. and Telgarsky, Matus and Yu, Bin},
  booktitle={Conference on Learning Theory},
  year={2024},
}

@inproceedings{zou2021benefits,
author={Difan Zou and Jingfeng Wu and Vladimir Braverman and Quanquan Gu and Dean P. Foster and Sham M. Kakade},
title={The Benefits of Implicit Regularization from {SGD} in Least Squares Problems},
booktitle={Neural Information Processing Systems},
year={2021}
}

@inproceedings{ji2020gradient,
  title={Gradient descent follows the regularization path for general losses},
  author={Ji, Ziwei and Dud{\'\i}k, Miroslav and Schapire, Robert E. and Telgarsky, Matus},
  booktitle={Conference on Learning Theory},
  year={2020}
}

@inproceedings{ji2019implicit,
  title={The implicit bias of gradient descent on nonseparable data},
  author={Ji, Ziwei and Telgarsky, Matus},
  booktitle={Conference on Learning Theory},
  year={2019}
}

@article{schapire1998boosting,
  title={Boosting the margin: {A} new explanation for the effectiveness of voting methods},
  author={Schapire, Robert E. and Freund, Yoav and Bartlett, Peter L. and Lee, Wee Sun},
  journal={Annals of Statistics},
  volume={26},
  number={5},
  pages={1651--1686},
  year={1998},
}

@article{rosset2004boosting,
  author = {Saharon Rosset and Ji Zhu and Trevor Hastie},
  journal = {Journal of Machine Learning Research},
  pages = {941--973},
  title = {Boosting as a regularized path to a maximum margin classifier},
  volume = 5,
  year = 2004}

@inproceedings{telgarsky2013margins,
  title={Margins, shrinkage, and boosting},
  author={Telgarsky, Matus},
  booktitle={International Conference on Machine Learning},
  year={2013},
}

@article{neyshabur2014search,
  title={In search of the real inductive bias: {On} the role of implicit regularization in deep learning},
  author={Neyshabur, Behnam and Tomioka, Ryota and Srebro, Nathan},
  journal={arXiv preprint arXiv:1412.6614},
  year={2014}
}

@article{soudry2018implicit,
  title={The implicit bias of gradient descent on separable data},
  author={Soudry, Daniel and Hoffer, Elad and Nacson, Mor Shpigel and Gunasekar, Suriya and Srebro, Nathan},
  journal={Journal of Machine Learning Research},
  volume={19},
  number={70},
  pages={1--57},
  year={2018}
}

@inproceedings{gunasekar2018characterizing,
  title={Characterizing implicit bias in terms of optimization geometry},
  author={Gunasekar, Suriya and Lee, Jason and Soudry, Daniel and Srebro, Nathan},
  booktitle={International Conference on Machine Learning},
  year={2018},
}

@article{sun2023unified,
  title={A unified approach to controlling implicit regularization via mirror descent},
  author={Sun, Haoyuan and Gatmiry, Khashayar and Ahn, Kwangjun and Azizan, Navid},
  journal={Journal of Machine Learning Research},
  volume={24},
  number={393},
  pages={1--58},
  year={2023}
}

@incollection{prechelt2002early,
  title={Early stopping-but when?},
  author={Prechelt, Lutz},
  booktitle={Neural Networks: {Tricks} of the trade},
  year={2002}
}

@article{buhlmann2003boosting,
  title={Boosting with the {L2} loss: {Regression} and classification},
  author={Buhlmann, Peter and Yu, Bin},
  journal={Journal of the American Statistical Association},
  volume={98},
  number={462},
  pages={324--339},
  year={2003},
}

@article{zhang2005boosting,
  author = {Tong Zhang and Bin Yu},
  journal = {Annals of Statistics},
  number = {4},
  pages = {1538--1579},
  title = {Boosting with early stopping: {Convergence} and consistency},
  volume = {33},
  year = {2005}}

@article{yao2007early,
  title={On early stopping in gradient descent learning},
  author={Yao, Yuan and Rosasco, Lorenzo and Caponnetto, Andrea},
  journal={Constructive Approximation},
  volume={26},
  number={2},
  pages={289--315},
  year={2007},
}

@article{raskutti2014early,
  title={Early stopping and non-parametric regression: {An} optimal data-dependent stopping rule},
  author={Raskutti, Garvesh and Wainwright, Martin J. and Yu, Bin},
  journal={Journal of Machine Learning Research},
  volume={15},
  number={11},
  pages={335--366},
  year={2014},
}

@inproceedings{wei2017early,
  title={Early stopping for kernel boosting algorithms: {A} general analysis with localized complexities},
  author={Wei, Yuting and Yang, Fanny and Wainwright, Martin J},
  booktitle={Neural Information Processing Systems},
  year={2017}
}

@inproceedings{suggala2018connecting,
  author = {Arun S. Suggala and Adarsh Prasad and Pradeep Ravikumar},
  booktitle = {Neural Information Processing Systems},
  title = {Connecting Optimization and Regularization Paths},
  year = {2018}}

@inproceedings{ali2019continuous,
  author = {Alnur Ali and J. Zico Kolter and Ryan J. Tibshirani},
  booktitle = {International Conference on Artificial Intelligence and Statistics},
  title = {A continuous-time view of early stopping for least squares},
  year = {2019}}

@inproceedings{ali2020implicit,
  title={The implicit regularization of stochastic gradient flow for least squares},
  author={Ali, Alnur and Dobriban, Edgar and Tibshirani, Ryan J.},
  booktitle={International Conference on Machine Learning},
  year={2020},
}

@article{sonthalia2024regularization,
title={On Regularization via Early Stopping for Least Squares},
author={Rishi Sonthalia and Jackie Lok and Elizaveta Rebrova}, 
journal={arXiv preprint arXiv:2406.04425},
year={2024}}

@article{lemaire1996asymptotical,
  title={An asymptotical variational principle associated with the steepest descent method for a convex function},
  author={Lemaire, Bastien},
  journal={Journal of Convex Analysis},
  volume={3},
  number={1},
  pages={63--70},
  year={1996},
}

@article{sason2015reverse,
  title={On reverse {Pinsker} inequalities},
  author={Sason, Igal},
  journal={arXiv preprint arXiv:1503.07118},
  year={2015}
}

@article{bauschke1997legendre,
  author = {Heinz H. Bauschke and Jonathan M. Borwein},
  journal = {Journal of Convex Analysis},
  number = 1,
  pages = {27--67},
  title = {Legendre Functions and the Method of Random {Bregman} Projections},
  volume = 4,
  year = 1997}

@book{vershynin2018high,
  title={High-Dimensional Probability: An Introduction with Applications in Data Science}, 
  author={Vershynin, Roman},
  year={2018},
  publisher={Cambridge}
}

@book{wainwright2019high,
  title={High-dimensional Statistics: A Non-Asymptotic Viewpoint},
  author={Wainwright, Martin J.},
  year={2019},
  publisher={Cambridge}
}

@book{wasserman2004all,
  author = {Larry Wasserman},
  title = {All of Statistics: A Concise Course in Statistical Inference},
  year = {2004},
  publisher = {Springer}}

@article{bach2010self,
  title={Self-concordant analysis for logistic regression},
  author={Bach, Francis},
  journal={Electronic Journal of Statistics},
  volume={4},
  pages={383--414},
  year={2010}
}

@article{ostrovskii2021finite,
  title={Finite-sample analysis of {M-estimators} using self-concordance},
  author={Ostrovskii, Dmitrii M. and Bach, Francis},
  journal={Electronic Journal of Statistics},
  volume={15},
  pages={326--391},
  year={2021}
}

@article{hsu2012tail,
  author = {Daniel Hsu and Sham Kakade and Tong Zhang},
  journal = {Electronic Communications in Probability},
  number = 52,
  pages = {1--6},
  title = {A tail inequality for quadratic forms of subgaussian random vectors},
  volume = 17,
  year = 2012}

@inproceedings{lin2017sharp,
  author = {Lin, Kevin and Sharpnack, James and Rinaldo, Alessandro and Tibshirani, Ryan J.},
  booktitle = {Neural Information Processing Systems},
  title = {A Sharp Error Analysis for the Fused Lasso, with Application to Approximate Changepoint Screening},
  year = 2017}

@article{pollard2017few,
  title={A few good inequalities},
  author={Pollard, David},
  year={2017},
  journal={Lecture Notes: Statistics 600: Advanced Probability},
  url={http://www.stat.yale.edu/~pollard/Courses/600.spring2017/Handouts/Basic.pdf}
}

@article{wolpert1992stacked,
  title={Stacked generalization},
  author={Wolpert, David H.},
  journal={Neural Networks},
  volume={5},
  number={2},
  pages={241--259},
  year={1992}
}

@article{breiman1996stacked,
  title={Stacked regressions},
  author={Breiman, Leo},
  journal={Machine Learning},
  volume={24},
  number={1},
  pages={49--64},
  year={1996}
}

@article{tibshirani2023empirical,
  title={Empirical process theory for nonparametric analysis},
  author={Tibshirani, Ryan J.},
  year={2023},
  journal={Lecture Notes, Stat 241B: Advanced Topics in Statistical Learning},
  url={https://www.stat.berkeley.edu/~ryantibs/statlearn-s23/lectures/emp_process.pdf},
}

@article{sidford2020convex,
  title={Smooth convex generalizations},
  author={Aaron Sidford},
  year={2020},
  journal={Lectures Notes, MSE 213: Introduction to Optimization Theory},
  url={https://web.stanford.edu/~sidford/courses/20fa_opt_theory/sidford_mse213_2020fa_chap_5_extensions.pdf},
}

@inproceedings{tsybakov2003optimal,
  title={Optimal rates of aggregation},
  author={Alexandre B. Tsybakov},
  booktitle={Conference on Learning Theory},
  year={2003},
}

@article{juditsky2005recursive,
  title={Recursive aggregation of estimators via the mirror descent algorithm with averaging},
  author={Anatoli B. Juditsky and Alexander V. Nazin and Alexandre B. Tsybakov and Nicolas Vayatis},
  journal={Methods of Signal Processing},
  volume={31},
  pages={368--384},
  year={2005},
}

@article{lecue2007optimal,
  title={Optimal rates of aggregation in classification under low noise assumption},
  author={Lecu{\'e}, Guillaume},
  journal={Bernoulli},
  volume={13},
  number={4},
  pages={1000--1022},
  year={2007},
}

@article{juditsky2008learning,
  title={Learning by mirror averaging},
  author={Juditsky, Anatoli B. and Rigollet, Philippe and Tsybakov, Alexandre B.},
  journal={Annals of Statistics},
  volume={36},
  number={5},
  pages={2183--2206},
  year={2008},
}

@inproceedings{dalalyan2009aggregation,
  title={Aggregation by exponential weighting and sharp oracle inequalities},
  author={Dalalyan, Arnak S. and Tsybakov, Alexandre B.},
  booktitle={Conference on Learning Theory},
  year={2009},
}

@article{lecue2013optimality,
  title={On the optimality of the aggregate with exponential weights for low temperatures},
  author={Lecu{\'e}, Guillaume and Mendelson, Shahar},
  journal={Bernoulli},
  volume={19},
  number={3},
  pages={646--675},
  year={2013},
}

@article{lecue2014optimal,
  title={Optimal learning with {Q-aggregation}},
  author={Lecu{\'e}, Guillaume and Rigollet, Philippe},
  journal={Annals of Statistics},
  volume={42},
  number={51},
  pages={211--224},
  year={2014},
}

@article{alquier2024user,
  title={User-friendly introduction to {PAC-Bayes} bounds},
  author={Alquier, Pierre},
  journal={Foundations and Trends in Machine Learning},
  volume={17},
  number={2},
  pages={174--303},
  year={2024},
}

@article{leung2006information,
  title={Information theory and mixing least-squares regressions},
  author={Leung, Gilbert and Barron, Andrew R.},
  journal={IEEE Transactions on Information Theory},
  volume={52},
  number={8},
  pages={3396--3410},
  year={2006},
}

@article{zhang2006information,
  title={Information-theoretic upper and lower bounds for statistical estimation},
  author={Zhang, Tong},
  journal={IEEE Transactions on Information Theory},
  volume={52},
  number={4},
  pages={1307--1321},
  year={2006},
}

@inproceedings{xu2017information,
  title={Information-theoretic analysis of generalization capability of learning algorithms},
  author={Xu, Aolin and Raginsky, Maxim},
  booktitle={Neural Information Processing Systems},
  year={2017},
}

@article{bauschke2017descent,
  title={A descent lemma beyond {Lipschitz} gradient continuity: {First-order} methods revisited and applications},
  author={Bauschke, Heinz H. and Bolte, J{\'e}r{\^o}me and Teboulle, Marc},
  journal={Mathematics of Operations Research},
  volume={42},
  number={2},
  pages={330--348},
  year={2017},
}

@article{efron2004least,
  author = {Bradley Efron and Trevor Hastie and Iain Johnstone and Robert Tibshirani},
  journal = {Annals of Statistics},
  number = 2,
  pages = {407--499},
  title = {Least angle regression},
  volume = 32,
  year = 2004}

\clearpage
\appendix
\renewcommand{\thesection}{\Alph{section}}
\renewcommand{\thetheorem}{\thesection\arabic{theorem}}
\renewcommand{\thelemma}{\thesection\arabic{lemma}}
\settocdepth{section}
\section{Definitions and notation}

We overview the standard definitions and notation used in the main paper. For a
set $S\subseteq\R^d$, we denote its interior and boundary by
\smash{$\interior(S)$} and \smash{$\partial S$}, respectively.  When $S$ is
finite, we denote its cardinality by $|S|$. 

For a convex set $\Omega$, a function is said to be convex if 
\smash{$f(\alpha x + (1-\alpha) y)\leq \alpha f(x) + (1-\alpha) f(y)$} for any 
\smash{$x,y\in\Omega$} and \smash{$\alpha\in[0,1]$}, and strictly convex if the 
inequality is strict for \smash{$x\neq y$} and \smash{$\alpha\in(0,1)$}. We
denote the subdifferential (set of subgradients) of $f$ at $x$ by
\smash{$\partial f(x)$}.  A function $f$ is said to be essentially strictly
convex if it is strictly convex on every convex subset of
\smash{$\{x \in \Omega : \partial f(x)\neq\emptyset\}$}. A function 
\smash{$f:\Omega\subseteq\R^d\to\R$} is said to be essentially smooth provided 
it satisfies three conditions: (i) \smash{$\interior(\Omega)\neq\emptyset$};
(ii) $f$ is differentiable on \smash{$\interior(\Omega)$}; and (iii)
\smash{$\lim_{k\to\infty}\|\nabla f(x_k)\|_2=\infty$} for any sequence
\smash{$\{x_k\}_{k=1}^\infty\subset\Omega$} converging to a point in
\smash{$\partial\Omega$}. Finally, a function $f$ is said to be of Legendre type
if it is both essentially smooth and essentially strictly convex.

A differentiable function $f$ is called $\alpha$-strongly convex with respect to
a norm \smash{$\|\cdot\|$}, for a given $\alpha>0$, if \smash{$f(y) \geq f(x) +
  \langle \nabla f(x), y-x\rangle +\frac\alpha2 \|x-y\|^2$}.  A differentiable
function $f$ is called $L$-smooth with respect to a norm \smash{$\|\cdot\|$},
for a given $L>0$, if the gradient map \smash{$\nabla f$} is $L$-Lipschitz with
respect to \smash{$\|\cdot\|$}, which recall means that {$\|\nabla f(x)-\nabla
  f(y)\|_* \leq L\|x-y\|$} for all \smash{$x,y$}, where \smash{$\|\cdot\|_*$}
denotes the dual norm.

\section{Proofs for Section \ref{sec:basic-ineq}}

\subsection{Descent lemmas}

For completeness, we state and prove the standard descent lemma for gradient
descent. 

\begin{lemma}
\label{lem:gd-descending}
Let \smash{$f:\R^d\to \R$} be convex, differentiable, and $L$-smooth for
$L>0$. For the gradient descent update in \eqref{eq:gd-iterate} with step size
\smash{$\eta_t \in(0,2/L]$}, 
\[
f(\theta_{t+1}) \leq f(\theta_t) - \eta_t\Big(1-\frac L2 \eta_t\Big) \|\nabla
f(\theta_t)\|_2^2 \leq f(\theta_t).
	\]
\end{lemma}

\begin{proof}
By $L$-smoothness,
\[
f(\theta_{t+1}) \leq f(\theta_t) + \langle \nabla f(\theta_t), 
\theta_{t+1}-\theta_t \rangle + \frac{L}{2} \|\theta_{t+1}-\theta_t\|_2^2. 
\]
Substituting \smash{$\theta_{t+1} = \theta_t - \eta_t \nabla f(\theta_t)$}
into the above and simplifying gives the desired result.
\end{proof}

Next we state and prove the corresponding descent lemma for mirror descent.

\begin{lemma}
\label{lem:md-descending}
Let $\Omega,\cC$ be closed, convex sets in $\R^d$ with nonempty interior, and 
$\cC\subseteq\Omega$. Assume $f:\Omega\to\R$ is convex on $\cC$, and
$L$-smooth with respect to a norm $\|\cdot\|$ on
$\cC\cap\interior(\Omega)$, for $L>0$. Assume $\phi:\Omega\to\R$ is of
Legendre type, and $\alpha$-strongly convex with respect to $\|\cdot\|$ on
$\Omega$, for $\alpha>0$. For the mirror descent update in \eqref{eq:md-iterate}
with step size $\eta_t\in(0,2\alpha/L]$, 
\[
f(\theta_{t+1}) \leq f(\theta_t) + \Big(\frac L2 - \frac\alpha{\eta_t}\Big)
\|\theta_t - \theta_{t+1}\|^2 \leq f(\theta_t).
\]
\end{lemma}

\begin{proof}
By $L$-smoothness,
\[
f(\theta_{t+1}) \leq f(\theta_t)  +  \langle \nabla f(\theta_t),
\theta_{t+1}-\theta_t \rangle + \frac{L}{2} \|\theta_{t+1}-\theta_t\|^2. 
\]
Meanwhile, by the three-point lemma and the first-order condition for convexity
(as derived in part 1 of the proof of Theorem \ref{thm:md-basicineq}), 
\[
\eta_t \langle \nabla f(\theta_t), \theta_{t+1} - z\rangle
\leq \Dphi(z,\theta_t) - \Dphi(z, \theta_{t+1}) - \Dphi(\theta_{t+1}, \theta_t).
\]
and plugging in $z = \theta_t$ gives 
\[
\eta_t \langle \nabla f(\theta_t), \theta_{t+1}- \theta_t\rangle 
\leq -\Dphi(\theta_t, \theta_{t+1}) - \Dphi(\theta_{t+1}, \theta_t).
\]
Combining this with the $L$-smoothness inequality, we obtain
\[
f(\theta_{t+1}) \leq f(\theta_t) -\frac{1}{\eta_t} \Dphi(\theta_t, \theta_{t+1}) 
- \frac{1}{\eta_t} \Dphi(\theta_{t+1}, \theta_t) + \frac L2
\|\theta_{t+1}-\theta_t\|^2. 
\]
By $\alpha$-strong convexity, the terms {$\Dphi(\theta_t, \theta_{t+1})$} and
{$\Dphi(\theta_{t+1}, \theta_{t})$} are each lower bounded by
\smash{$\frac\alpha2\|\theta_{t+1}-\theta_t\|^2$}. Plugging this in and
simplifying gives the desired result.  
\end{proof}

\subsection{Proof of Corollary \ref{cor:egd-envelope}} 

We will establish the two upper bounds given in \eqref{eq:egd-envelope} 
separately. For the first upper bound, which involves the penalty 
\smash{$\frac{d+1}{2}\|\pi-z\|_1^2$}, this is a consequence of
\eqref{eq:md-envelope} after checking the appropriate strong convexity and
smoothness properties of the negative entropy $\phi$. By Pinsker's inequality,
we know that $\alpha=1$ is the strong convexity parameter of $\phi$ with respect
to $\|\cdot\|_1$. For smoothness, we can calculate, letting $s = z-\pi$,       
\begin{align*}
\DKL(z,\pi) &= \sum_{i=1}^d z_i \log(dz_i) \\
&= \sum_{i=1}^d (s_i + 1/d)  \log(1+d s_i) \\
&\leq \sum_{i=1}^d (s_i + 1/d) d s_i \\
&= d\sum_{i=1}^d s_i^2,
\end{align*}
where the inequality uses \smash{$\log(1+x) \leq x$} for $x \geq -1$. Meanwhile,
denoting \smash{$C=\sum_{i=1}^d s_i \cdot \one (s_i\geq 0)$}, note that
\smash{$\sum_{i=1}^d s_i\cdot \one (s_i< 0) = -C$} and $\|s\|_1=2C$ since
\smash{$\sum_{i=1}^d s_i=0$}; therefore  
\begin{align*}
\sum_{i=1}^d s_i^2 
&= \sum_{i=1}^d s_i^2 \cdot \one (s_i\geq 0) + 
\sum_{i=1}^d s_i^2 \cdot \one (s_i< 0) \\
&\leq \bigg(\sum_{i=1}^d |s_i| \cdot \one (s_i\geq 0)\bigg)^2 +
\bigg(\sum_{i=1}^d |s_i| \cdot \one (s_i< 0)\bigg)^2 \\
&= 2C^2 \\
&= \|s\|_1^2/2.
\end{align*}
Combining this with the second-to-last display therefore gives 
\[
\DKL(z,\pi) \leq \frac{d}{2} \|s\|_1^2,
\]
which proves the desired smoothness result with $G=d$. 

For the second upper bound, which involves the penalty \smash{$\frac{1}{2} \|\pi
  - z\|_1^2 + \frac{\log d}{2} \|\pi - z\|_1$}, we now use a reverse
Pinsker-type inequality due to Theorem 1 in \citet{sason2015reverse}:
\[
\DKL(z,\pi)\leq\frac{\log d}{2(1-1/d)}\|z-\pi\|_1.
\] 
By strong convexity, the triangle inequality, and Young's inequality,  
\[
\DKL(z, \theta_T) \geq \frac14 \|\pi-\theta_T\|_1^2 - \frac12 \|z-\pi\|_1^2.
\]
Combining the previous two bounds gives 
\[
\DKL(z,\pi)-\DKL(z,\theta_T) \leq \frac{\log d}{2(1-1/d)}\|z-\pi\|_1 
- \frac14 \|\pi-\theta_T\|_1^2 + \frac12 \|z-\pi\|_1^2.
\]
Applying \eqref{eq:md-basicineq} and rearranging proves the desired result. 

\subsection{Proof of Corollary \ref{cor:md-optlemaire}}

Parts 1 and 2 follow from arguments analogous to those for Corollary 
\ref{cor:gd-optlemaire}. We prove parts 3 and 4 below.

\medskip\noindent
\underline{Part 3.}
By part 2 \smash{$\{\Dphi(s,\theta_t)\}_{t=0}^\infty$} is nonincreasing, and by 
strong convexity \smash{$\|s-\theta_t\|^2 \leq \frac2\alpha
  D_\phi(s,\theta_t)$} for each $t$. Thus by the Bolzano-Weierstrass theorem
there exists a subsequence which converges, i.e., \smash{$\lim_{k\to\infty}
  \theta_{t_k} = \theta_\infty$}. By part 1 and continuity of $f$, we know that
$\theta_\infty \in S$.  

To see that the whole sequence \smash{$\{\theta_t\}_{t=0}^\infty$} converges,
suppose not. Then there exists $\delta>0$ and a subsequence
\smash{$\{\theta_{t_k}\}_{k=0}^\infty$} such that
\smash{$\|\theta_{t_k}-\theta_\infty\| \geq \delta$} for all $k$. With the same
argument as used above, we know there exists a subsubsequence
\smash{$\{\theta_{t_{k_m}}\}_{m=0}^\infty$} which converges to another limit
point \smash{$\tilde\theta_\infty \in S$} such that
\smash{$\|\tilde\theta_\infty - \theta_\infty\| \geq \delta$}. Now let case (i)
denote the assumption that $S \cap \interior(\Omega) \neq \emptyset$, and case
(ii) denote the assumption that $\Dphi$ is continuous in its second argument
relative to $\interior(\Omega)$. Consider the following. 
\begin{itemize}
\item In case (i), picking any $s \in S \cap \interior(\Omega)$, since
  \smash{$\{\Dphi(s,\theta_{t_k})\}_{k=0}^\infty$} is nonincreasing and hence 
  bounded, we know by Theorem 3.8(ii) in \citet{bauschke1997legendre} that
  $\theta_\infty \in S$ and \smash{$\Dphi(\theta_\infty,\theta_{t_k}) \to 0$}.
\item In case (ii), by definition of continuity, and since each
  \smash{$\theta_{t_k} \in \interior(\Omega)$} (due to the nature of the 
  mirror descent updates and the Legendre property of $\phi$), we have 
  \smash{$\Dphi(\theta_\infty,\theta_{t_k}) \to 0$}.  
\end{itemize}
In either case, \smash{$\Dphi(\theta_\infty,\theta_{t_k}) \to 0$} as 
$k \to \infty$. By the nonincreasing property of
\smash{$\{\Dphi(s,\theta_t)\}_{t=0}^\infty$}, we have that
\smash{$\Dphi(\theta_\infty,\theta_t) \to 0$} as $t \to \infty$. By a similar
argument applied to the subsubsequence
\smash{$\{\theta_{t_{k_m}}\}_{m=0}^\infty$}, we hence also
\smash{$\Dphi(\tilde\theta_\infty,\theta_t) \to 0$} as $t \to \infty$. But this
leads to the desired contradiction, as the same sequence cannot converge to  
distinct points \smash{$\theta_\infty,\tilde\theta_\infty$}. Formally, 
\[
\max\Big\{\Dphi(\theta_\infty, \theta_t), \Dphi(\tilde\theta_\infty,
\theta_t)\Big\} \geq \max \Big\{\frac\alpha2 \|\theta_\infty - \theta_t\|^2,
\frac\alpha2 \|\tilde\theta_\infty - \theta_t\|^2\Big\}	\geq \frac\alpha2
\Big(\frac\delta2\Big)^2. 
\]
Hence, the entire sequence \smash{$\{\theta_t\}_{t=0}^\infty$} converges,
and \smash{$\lim_{t\to\infty}\theta_t = \theta_\infty \in S$}.

\medskip\noindent
\underline{Part 4.}
Let \smash{$p = \proj^\phi_S(\theta_0)$} and $v = p - \theta_\infty$. Assume $v
\neq 0$. For any $c \in \R$, note that $p+cv \in S$ since $S$ affine. 
We can obtain three relations regarding $p$, \smash{$\theta_\infty$}, and
\smash{$p+cv$}. First, the three-point identity for $\Dphi$ gives 
\[
\Dphi(p + cv, \theta_\infty) - \Dphi(p + cv, p) - \Dphi(p, \theta_\infty) 
= \langle\nabla\phi(p)-\nabla\phi(\theta_\infty), p+ cv - p\rangle.
\]
Second, by part 2, we know \smash{$\Dphi(p + c v, \theta_\infty) \leq \Dphi(p +
  cv, \theta_0)$}. Third, as $S$ is affine, the generalized Pythagorean theorem 
for Bregman projection holds with equality:
\[
\Dphi(p + cv, \theta_0) = \Dphi(p + cv, p) + \Dphi(p, \theta_0). 
\]
Combining these three relations, we have
\begin{align*}
\langle\nabla\phi(p)-\nabla\phi(\theta_\infty), c v \rangle
&= \Dphi(p+cv, \theta_\infty) - \Dphi(p+cv, p) - \Dphi(p, \theta_\infty) \\ 
&\leq \Dphi(p+cv, \theta_0) - \Dphi(p+cv, p) - \Dphi(p, \theta_\infty) \\  
&= \Dphi(p,\theta_0) - \Dphi(p,\theta_\infty).
\end{align*}
As this holds for all $c \in \R$, we conclude
\smash{$\langle\nabla\phi(p)-\nabla\phi(\theta_\infty), v\rangle=0$}. By
$\alpha$-strong convexity, 
\[
0 = \langle\nabla\phi(p)-\nabla\phi(\theta_\infty), p - \theta_\infty \rangle 
\geq \alpha\|p-\theta_\infty\|^2, 
\]
so $v=0$ and $p = \theta_\infty$. This completes the proof.

\section{Proofs for Section \ref{sec:glm}}

\allowdisplaybreaks
\subsection{Proof of Proposition \ref{prop:risk-ridgeglm}}

From \eqref{eq:glm-risk-diff}, we can apply Young's inequality \smash{$2ab \leq
  ca^2+b^2/c$}, with \smash{$c=\lambda$}, to the first term on the right-hand
side, yielding  
\begin{align*}
\frac1n \epsilon^\top X(\hthetalambda -\theta) + \lambda \Big(\|\theta\|_2^2 -
  \|\hthetalambda\|_2^2 \Big) 
&\leq \frac{1}{2\lambda} \Big\|\frac{X^\top\epsilon}{n}\Big\|_2^2
  + \frac{\lambda}{2} \|\hthetalambda-\theta\|_2^2 + 
  \lambda \Big(\|\theta\|_2^2 - \|\hthetalambda\|_2^2 \Big) \\ 
&\leq \frac{1}{2\lambda} \Big\|\frac{X^\top\epsilon}{n}\Big\|_2^2
  + \lambda\Big(\|\hthetalambda\|_2^2 + \|\theta\|_2^2 \Big)
  + \lambda \Big(\|\theta\|_2^2 - \|\hthetalambda\|_2^2 \Big) \\
&= \frac{1}{2\lambda} \Big\|\frac{X^\top\epsilon}{n}\Big\|_2^2
  + 2\lambda \|\theta\|_2^2.
\end{align*}
This proves the first claim in the proposition.

To prove the second claim, since the $\epsilon_i$ are independent and each
$\epsilon_i\sim\sG(\sigma^2)$, Theorem 1 of \citet{hsu2012tail} implies 
\begin{equation} 
\label{eq:hsu-lemma}
\P\bigg(\Big\| \frac{X^\top \epsilon}{n} \Big\|_2^2 > \frac{\sigma^2}{n}
\Big[ \trhSigma + 2 \sqrt{\delta} \FhSigma + 2 \delta \ophSigma \Big]\bigg)    
\leq e^{-\delta}. 
\end{equation}
Applying this to the first result of the proposition gives 
\[
\risk(\hthetalambda) - \min_{\theta: \|\theta\|_2\leq b} \risk(\theta)
\leq \frac{R}{2\lambda} + 2\lambda b^2,
\]
with probability at least \smash{$1-e^{-\delta}$}, where
\begin{equation} 
\label{eq:R-def}
R = \frac{\sigma^2}{n} \Big[ \trhSigma + 2 \sqrt{\delta} \FhSigma + 2
\delta \ophSigma \Big]. 
\end{equation}
Choosing \smash{$\lambda = \sqrt{R} / (2b)$}, we get that
\[
\risk(\hthetalambda) - \min_{\theta: \|\theta\|_2\leq b} \risk(\theta)
\leq 2 b \sqrt{R},
\]
with probability at least \smash{$1-e^{-\delta}$}, which completes the proof. 
	
\subsection{Proof of Theorem \ref{thm:glm-ridge-risk}} 

In the Gaussian case \smash{$\epsilon_i\sim \sG(\sigma_i^2)$}, and in the
Bernoulli case \smash{$\epsilon_i \sim \sG(1/4)$}. The result in either of
these cases is then a direct application of Proposition
\ref{prop:risk-ridgeglm}. In the Poisson case, the noise $\epsilon_i =
y_i-\mu_i$ is not sub-Gaussian so we will use a truncation argument, similar to
that in Appendix A.4 of \citet{lin2017sharp}. Define an event 
\[
\cE = \{\epsilon_i < M, \; i=1,\dots,n\},
\quad \text{where} \quad  
M = 4\big( \|\boldmu\|_\infty + 1/3\big) \log n.
\]
Note that \smash{$M>1$} for \smash{$n\geq 3$}. Further,
\[
\P(\cE^c) \leq \sum_{i=1}^n \P(y_i - \mu_{i} \geq M).
\]
Now using a Poisson concentration result from \citet{pollard2017few}, for
\smash{$X\sim\text{Pois}(\mu)$}, and for all $x>0$, 
\[
\P( X-\mu \geq x) \leq \exp \bigg(-\frac{x^2}{2\mu}
\psi\Big(\frac{x}{\mu}\Big)\bigg),
\quad \text{where} \quad  
\psi(x) = \frac{(1+x)\log(1+x)-x}{x^2/2}.
\]
Moreover, when \smash{$x\geq 1$},
\[
\frac{x^2}{2\mu} \psi\Big(\frac{x}{\mu}\Big) \geq \frac{1/2}{\mu+1/3} x.  
\]
Therefore, we have the following for each $i$, 
\[
\P(\epsilon_i \geq M) \leq \exp\bigg( -\frac{1/2}{\mu_{i} + 1/3} M \bigg)   
\leq \exp \bigg( -\frac{1/2}{\|\boldmu\|_{\infty}+ 1/3} M \bigg) =  1/n^2,
\]
which means that \smash{$\P(\cE^c) \leq \sum_{i=1}^n 1/n^2 = 1/n$}. Now define
an event 
\[
\cS_\delta = \bigg\{ \Big\|\frac{X^\top \epsilon}{n} \Big\|_2^2 >
\frac{\sigma^2}{n} \Big[\trhSigma + 2 \sqrt{\delta} \FhSigma + 2 \delta
\ophSigma \Big] \bigg\},
\]
where $\sigma^2$ is yet to be specified. Note $\P(\cS_\delta) = \P(\cS_\delta
\cap \cE^c) + \P(\cS_\delta \cap \cE) \leq \P(\cE^c) + \P(\cS_\delta \cap \cE)
\leq  1/n + \P(\cS_\delta \cap \cE)$. It remains to bound the latter
probability. 

Define \smash{$\tilde\epsilon_i = \epsilon_i \one(\epsilon_i < M)$}, and
\smash{$b_i = \E[\tilde\epsilon_i]$}. As $\E[\epsilon_i]=0$, we have
\smash{$b_i = -\E[\epsilon_i \one(\epsilon_i \geq M)]$}. By Cauchy-Schwarz,
\smash{$|b_i|^2 \leq \mu_i \P(\epsilon_i \geq M) \leq \mu_i/n^2$}. Finally, let
\smash{$\bar\epsilon_i = \tilde\epsilon_i - b_i$}, and \smash{$\bar{M} = M +
  \sqrt{\|\mu\|_\infty}/n$}. Then \smash{$\bar\epsilon_i$} has mean zero and is 
bounded in magnitude by \smash{$\bar{M}$}, and so it is sub-Gaussian with
parameter \smash{$\bar\sigma^2 = \bar{M}^2$}. 

This means that we can apply \eqref{eq:hsu-lemma} to bound the tail probability
of \smash{$\|X^\top\bar\epsilon/n\|_2^2$}, concretely
\[
\P\bigg(\Big\| \frac{X^\top\bar\epsilon}{n} \Big\|_2^2 > \frac{\bar\sigma^2}{n} 
\Big[ \trhSigma + 2 \sqrt{\delta} \FhSigma + 2 \delta \ophSigma \Big]\bigg)    
\leq e^{-\delta}.
\]
Now observe
\begin{align*}
\Big\| \frac{X^\top\tilde\epsilon}{n} \Big\|_2^2
&\leq 2\Big\| \frac{X^\top\bar\epsilon}{n} \Big\|_2^2 +
2\Big\| \frac{X^\top b}{n} \Big\|_2^2 \\ 
&\leq 2\Big\| \frac{X^\top\bar\epsilon}{n} \Big\|_2^2 +
\frac{2}{n} \|b\|_2^2 \ophSigma \\
&\leq 2\Big\| \frac{X^\top\bar\epsilon}{n} \Big\|_2^2 +
\frac{2}{n^2} \|\mu\|_\infty \ophSigma.
\end{align*}
This means that 
\[
\P\bigg(\Big\| \frac{X^\top\tilde\epsilon}{n} \Big\|_2^2 >
\frac{2\bar\sigma^2}{n} \Big[ \trhSigma + 2 \sqrt{\delta} \FhSigma + 2 \delta
\ophSigma \Big] + \frac{2\|\mu\|_\infty}{n^2} \ophSigma \bigg) \leq
e^{-\delta}. 
\]
Provided $\delta \geq 1$,
\begin{align*}
\frac{2\bar\sigma^2}{n} \Big[ \trhSigma + 2 \sqrt{\delta} \FhSigma + 2 \delta
\ophSigma \Big] + \frac{2\|\mu\|_\infty}{n^2} \ophSigma 
&\leq \frac{2\bar\sigma^2+\|\mu\|_\infty/n}{n} \Big[ \trhSigma + 2 \sqrt{\delta}
\FhSigma + 2 \delta \ophSigma \Big] \\
&\leq \frac{\sigma^2}{n} \Big[ \trhSigma + 2 \sqrt{\delta} \FhSigma + 2 \delta
\ophSigma \Big],
\end{align*}
where the last line follows by calculating and defining
\smash{$\sqrt{2\bar\sigma^2+\|\mu\|_\infty/n} \leq
  \sqrt{2}(M+\sqrt{\|\mu\|_\infty}/n)+\sqrt{\|\mu\|_\infty/n}$} 
\smash{$\leq 6(\|\mu\|_\infty+1/3)\log{n} + 3 \sqrt{\|\mu\|_\infty} 
  := \sigma$}. Putting these pieces together, we have shown  
\[
\P\bigg(\Big\| \frac{X^\top\tilde\epsilon}{n} \Big\|_2^2 >
\frac{\sigma^2}{n} \Big[ \trhSigma + 2 \sqrt{\delta} \FhSigma + 2 \delta
\ophSigma \Big] \bigg) \leq e^{-\delta}.
\]
Note that on $\cE$, each \smash{$\tilde\epsilon_i = \epsilon_i$}, and therefore
$\P(\cS_\delta \cap \cE) \leq e^{-\delta}$ as well. This proves that
$\P(\cS_\delta) \leq 1/n + e^{-\delta}$, as desired. The remainder of the proof
follows by choosing $\lambda$ as before, and bounding the subsequent excess 
risk. 

\subsection{Proof of Proposition \ref{prop:risk-gdglm}}

From \eqref{eq:glm-risk-diff-gd}, we can apply Young's inequality \smash{$2ab
  \leq ca^2+b^2/c$}, with \smash{$c=\lambda$}, to the first term on the
right-hand side. This proves the first claim in the proposition.  

To prove the second claim, combining \eqref{eq:hsu-lemma} with the first result
of the proposition gives 
\[
\risk(\thetagdT) - \min_{\theta: \|\theta\|_2\leq b} \risk (\theta) \leq   
\frac{1}{2\lambda_T} R +\frac{\lambda_T}{2} b^2,
\]
with probability at least \smash{$1-e^{-\delta}$}, where \smash{$R$} is defined
in \eqref{eq:R-def}. At \smash{$\lambda^*_T = \sqrt{R}/b$}, this becomes
\[
\risk(\thetagdT) - \min_{\theta: \|\theta\|_2\leq b} \risk (\theta) \leq
b \sqrt{R}, 
\]
with probability at least \smash{$1-e^{-\delta}$}. In the case
\smash{$\lambda^*_T = \sqrt{R}/b$} is obtainable by $\lambda_T = 1/(\eta T)$ for
an integer $T$, this establishes the result in \eqref{eq:glm-gd-riskhighprob}
for gradient descent. Otherwise, we set \smash{$T = \lceil \frac{1}{\eta
  \lambda^*_T} \rceil$}; then \smash{$1/(\eta\lambda_T) \leq
  1/(\eta\lambda^*_T) + 1$}, i.e., \smash{$1/\lambda_T \leq 1/\lambda^*_T +
  \eta$}, so by Lemma \ref{lem:aux-gd-delta} (presented below),    
\[
\Big(\frac{1}{2\lambda_T} R + \frac{\lambda_T}{2} b^2\Big) -
\Big(\frac{1}{2\lambda^*_T} R + \frac{\lambda^*_T}{2} b^2\Big) \leq  
\frac{\eta R}{2}. 
\]
Thus in this case, we have an additional discretization error term of at most
\smash{$e_T = \eta R / 2$}. This completes the proof of the proposition.  

\begin{lemma} 
\label{lem:aux-gd-delta}
Consider a function \smash{$g(x) = \frac{a}{x} + bx$} on \smash{$(0,\infty)$}
with \smash{$a,b>0$}. This is minimized at \smash{$x^*= \sqrt{a/b}$}. If
\smash{$y>0$} satisfies \smash{$1/y = 1/x^* + c$} with \smash{$c\geq0$}, then 
\[
g(y) - g(x^*) = a c^2 y \leq ac.
\]
\end{lemma}

\begin{proof}
Observe
\begin{align*}
g(y) - g(x^*) 
&= a \Big(\frac1y - \frac{1}{x^*}\Big) + b(y-x^*) \\
&= ac - b \Big(\frac1y - \frac{1}{x^*}\Big)  yx^* \\
&= c (a-byx^*) \\
&= c \Big(a - \frac{ay}{x^*}\Big) \\
&= ac^2y \\
&\leq ac,
\end{align*}
where the fourth line uses \smash{$bx^* = a / x^*$}, and the last uses
\smash{$cy =(1/y - 1/x^*) y \leq 1$}.  
\end{proof}

\subsection{Proof of Theorem \ref{thm:glm-gd-risk}}

The theorem follows from Proposition \ref{prop:risk-gdglm} in the same way that 
Theorem \ref{thm:glm-ridge-risk} follows from Proposition
\ref{prop:risk-ridgeglm}. It only remains to identify (local) Lipschitz
parameters $L$ for the individual GLMs. Note that in general 
\[
\nabla^2 f(\theta) = \frac1n X^\top \nabla^2\bA(X\theta) X
= \frac1n \sum_{i=1}^n x_i A''(x_i^\top \theta) x_i^\top.
\]
\begin{itemize}
\item 
For case 1: Gaussian, we have \smash{$A(u) = u^2/2$} and \smash{$A''(u) =
  1$}. Hence \smash{$\nabla^2 f(\theta) = \hSigma$} and \smash{$L = \|
  \hSigma \|_\op$} is the global Lipschitz parameter. 
\item
For case 2: Bernoulli, we have \smash{$A(u) = \log(1+e^u)$} and \smash{$A''(u) =
  e^u/(1+e^u)^2 \leq \frac14$}. Hence \smash{$\nabla^2 f(\theta) \preceq
  \frac14 \hSigma$} and \smash{$L = \frac14 \| \hSigma \|_\op$} is the  
  global Lipschitz parameter. 
\item
For case 3: Poisson, we have \smash{$A(u)=e^u$} and \smash{$A''(u)=e^u$}. 
Note that \smash{$x_i^\top \theta \leq \|x_i\|_2 \|\theta\|_2 = bm$}, for
$\|\theta\|_2 \leq b$ and \smash{$m = \max_{i=1,\dots,n} \|x_i\|_2$}. Hence
\smash{$\nabla^2 f(\theta) \preceq e^{bm} \hSigma$} and \smash{$L = e^{bm}
  \|\hSigma\|_\op$} is the local Lipschitz parameter on \smash{$\ball_d(b)$}. 
\end{itemize}

\section{Proofs for Section \ref{sec:model-agg}}

\subsection{Proof of Proposition \ref{prop:risk-klglm}}

The first claim is a special case of a more general result for
Bregman-divergence-regularized GLMs, given below in Proposition 
\ref{prop:risk-bregglm}. In particular, this claim follows by setting $\phi$
to be the negative entropy function on $\Omega=\cC=\simplex$. This is 
1-strongly convex with respect to  $\|\cdot\|_1$ by Pinsker's inequality.    

To prove the second claim, note that each $X_j^\top \epsilon \sim
\sG(\sigma^2\|X_j\|_2^2)$. Thus, by a standard maximal inequality for 
sub-Gaussian random variables \citep{tibshirani2023empirical}, 
\begin{equation}
\label{eq:sg-maximal}
\P\Big( \| X^\top\epsilon \|_\infty > \sigma \sqrt{2n(\log(2d)+\delta)} \Big)
\leq e^{-\delta}.  
\end{equation}
Applying this to the first result of the proposition gives 
\[
\risk(\hthetalambda) - \min_{\theta: \DKL(\theta,u) \leq b} \risk(\theta) 
\leq \frac{R^2}{\lambda n^2} + 2\lambda b,
\]
with probability at least \smash{$1-e^{-\delta}$}, where
\begin{equation} 
\label{eq:R2-def}
R = \sigma \sqrt{2n(\log(2d)+\delta)}.
\end{equation}
Choosing \smash{$\lambda = R/(\sqrt{2b} n)$}, we get that
\[
\risk(\hthetalambda) - \min_{\theta: \DKL(\theta,u) \leq b} \risk(\theta)
\leq \frac{2\sqrt{2b} R}{n},
\]
with probability at least \smash{$1-e^{-\delta}$}, which completes the proof.  

\begin{proposition}
\label{prop:risk-bregglm}
Assume that $\phi : \Omega \to \R$ is $\alpha$-strongly convex on $\cC$
with respect to a norm \smash{$\| \cdot \|$}, for \smash{$\alpha > 0$},
and let \smash{$\|\cdot\|_*$} be the corresponding dual norm. Consider the
regularized GLM estimator defined by augmenting the GLM loss
\eqref{eq:glm-loss} with a Bregman divergence penalty, denoted 
\[
\hthetalambda =  \argmin_{\theta \in \cC} \Big( f(\theta) + \lambda
\Dphi(\theta,u) \Big),
\]
for \smash{$\lambda>0$}, and an anchor point \smash{$u\in\cC$}. Then for 
any reference point \smash{$\theta\in\cC$}, the prediction risk of
\smash{$\hthetalambda$} satisfies
\[
\risk(\hthetalambda)	\leq \risk(\theta)	+ \frac{1}{\alpha\lambda}
\Big\|\frac{X^\top\epsilon}{n}\Big\|_*^2 + 2\lambda \Dphi(\theta,u). 
\]
\end{proposition}

\begin{proof}
Following analogous steps to those leading up to \eqref{eq:glm-risk-diff},  
\[
\risk(\hthetalambda) - \risk(\theta)	\leq \frac1n \epsilon^\top X
(\hthetalambda-\theta) + \lambda \Big( \Dphi(\theta,u) - \Dphi(\hthetalambda,u)
\Big). 
\]
Using H{\"o}lder's inequality (with respect to the pair $\|\cdot\|$ and
$\|\cdot\|_*$) and Young's inequality $2ab \leq ca^2+b^2/c$, with
\smash{$c=\alpha\lambda$}, we have   
\begin{align*}
\frac1n \epsilon^\top X (\hthetalambda-\theta)
&\leq \Big\|\frac{X^\top\epsilon}{n}\Big\|_* \|\hthetalambda-\theta\| \\
&\leq \Big\|\frac{X^\top\epsilon}{n}\Big\|_* \|\hthetalambda-u\|
  + \Big\|\frac{X^\top\epsilon}{n}\Big\|_* \|\theta-u\| \\
&\leq \frac{1}{2 \alpha\lambda}\Big\|\frac{X^\top\epsilon}{n}\Big\|_*^2 
  + \frac{\alpha\lambda}{2} \|\hthetalambda-u\|^2
  + \Big\|\frac{X^\top\epsilon}{n}\Big\|_* \|\theta-u\|.
\end{align*} 
By strong convexity \smash{$\Dphi(u,v)\geq (\alpha/2) \|u-v\|^2$}, we obtain 
\begin{align*}
\risk(\hthetalambda)-\risk(\theta)
&\leq \frac{1}{2\alpha\lambda}\Big\|\frac{X^\top\epsilon}{n}\Big\|_*^2
  + \frac{\alpha\lambda}{2}\|\hthetalambda-u\|^2 
  + \Big\|\frac{X^\top\epsilon}{n}\Big\|_* \|\theta-u\| 
  + \lambda\Big(\Dphi(\theta,u) - \Dphi(\hthetalambda,u)\Big) \\ 
&\leq \frac{1}{2\alpha\lambda}\Big\|\frac{X^\top\epsilon}{n}\Big\|_*^2
  + \Big\|\frac{X^\top\epsilon}{n}\Big\|_* \|\theta-u\| 
  + \lambda \Dphi(\theta,u) \\
&\leq \frac{1}{\alpha\lambda}\Big\|\frac{X^\top\epsilon}{n}\Big\|_*^2
  + \frac{\alpha\lambda}{2}\|\theta-u\|^2 
  + \lambda \Dphi(\theta,u) \\
&= \frac{1}{\alpha\lambda}\Big\|\frac{X^\top\epsilon}{n}\Big\|_*^2
  + 2 \lambda \Dphi(\theta,u),
\end{align*}
where the second-to-last step uses Young's inequality again, and the last step 
uses strong convexity. 
\end{proof}

\subsection{Proof of Theorem \ref{thm:glm-kl-risk}}

This follows precisely as in the proof of Theorem \ref{thm:glm-ridge-risk}: for
the Gaussian and Bernoulli cases we simply identify $\sigma$ in the sub-Gaussian 
parameterization, and for the Poisson case, we use a truncation argument. 

\subsection{Proof of Proposition \ref{prop:risk-egdglm}}

The first claim is again a special case of a more general result for mirror
descent on a GLM loss, given below in Proposition \ref{prop:risk-mdglm}, applied
to the case where $\phi$ is the negative entropy function on 
$\Omega=\cC=\simplex$. The second claim follows from arguments just as in the
proof of Proposition \ref{prop:risk-klglm}, using the sub-Gaussian maximal
inequality \eqref{eq:sg-maximal}. This gives 
\[
\risk(\thetaegdT) - \min_{\theta: \DKL(\theta,u) \leq b} \risk (\theta) \leq
\frac{R^2}{2\lambda n^2} + \lambda b,
\]
with probability at least \smash{$1-e^{-\delta}$}, where \smash{$R$} is defined
in \eqref{eq:R2-def}. At \smash{$\lambda^*_T = R/(\sqrt{2b} n)$}, this gives 
\[
\risk(\thetaegdT) - \min_{\theta: \DKL(\theta,u) \leq b} \risk(\theta)
\leq \frac{\sqrt{2b} R}{n},
\]
which proves the result in \eqref{eq:glm-egd-riskhighprob} for the case of
integral $T^*$. For non-integral $T^*$, we then apply an analogous
discretization argument based on Lemma \ref{lem:aux-gd-delta}, which gives  
\[
\Big(\frac{R^2}{2\lambda_T n^2} + \lambda_T b\Big) - 
\Big(\frac{R^2}{2\lambda_T^* n^2} + \lambda_T^* b\Big) \leq \frac{R^2}{2n^2}
\eta.
\]
This verifies the additional discretization error term of at most \smash{$e_T =
  \eta R^2 / (2n^2)$}, and completes the proof.  

\begin{proposition}
\label{prop:risk-mdglm}
Under the assumptions of Proposition \ref{prop:risk-bregglm}, assume
additionally also that the GLM loss $f$ is $L$-smooth on $\cC$. Consider mirror descent
\eqref{eq:md-iterate} with $\eta_t = \eta \in (0,1/L]$, $t=0,1,2,\dots$, 
initialized at $\theta_0=u\in\cC$. Then for any reference point $\theta\in\cC$,
and for any $T \geq 1$ and $\lambda_T = 1/(\eta T)$,  
\[
\risk (\thetamdT) - \risk (\theta) \leq \frac{1}{2 \alpha \lambda_T}
\Big\|\frac{X^\top\epsilon}{n} \Big\|_*^2 + \lambda_T \Dphi(\theta,u). 
\]
\end{proposition}

\begin{proof}
By the basic inequality \eqref{eq:md-basicineq} in Theorem
\ref{thm:md-basicineq}, 
\[
f (\thetamdT) - f (\theta) \leq \lambda_T \Dphi (\theta,u) - \lambda_T 
\Dphi (\theta, \thetamdT).  
\]
Following analogous steps to those leading up to \eqref{eq:glm-risk-diff-gd},
and then arguments as in the proof of Proposition \ref{prop:risk-bregglm},
\begin{align*}
\risk (\thetamdT) - \risk (\theta)
& \leq \frac1n \epsilon^\top X (\thetamdT - \theta)
  + \lambda_T \Big( \Dphi (\theta,u) - \Dphi (\theta, \thetamdT) \Big) \\
&\leq \Big\| \frac{X^\top \epsilon}{n} \Big\|_* \|\thetamdT - \theta\|
  + \lambda_T \Big(\Dphi (\theta,u) - \Dphi (\theta, \thetamdT)  \Big) \\
&\leq \frac{1}{2\alpha\lambda_T} \Big\| \frac{X^\top \epsilon}{n} \Big\|_*^2
  + \frac{\alpha\lambda_T}{2} \|\thetamdT - \theta\|^2 + 
  \lambda_T \Big( \Dphi (\theta,u) - \Dphi (\theta, \thetamdT) \Big) \\
&\leq \frac{1}{2\alpha\lambda_T } \Big\| \frac{X^\top \epsilon}{n} \Big\|_*^2 
  + \lambda_T \Dphi (\theta,u).
\end{align*}
\end{proof}

\subsection{Proof of Theorem \ref{thm:glm-egd-risk}} 

After identifying the Lipschitz constant for each GLM loss, the arguments are
the identical as those in the proof of Theorem \ref{thm:glm-gd-risk}. Recall the
following equivalent characterization of $L$-smoothness with respect to
$\|\cdot\|_1$ \citep{sidford2020convex}:
\[
|v^\top \nabla^2 f(\theta) v| \leq L \|v\|_1^2,
\]
for all $\theta\in\simplex$ and directions $v\in\R^d$. Furthermore, 
\[
|v^\top \nabla^2 f(\theta) v| = \bigg| \frac1n \sum_{i=1}^n A''(x_i^\top \theta)  
(x_i^\top v)^2 \bigg| \leq \frac{\|v\|_1^2}{n} \sum_{i=1}^n A''(x_i^\top \theta)
\|x_i\|_\infty^2, 
\]
by H{\"o}lder's inequality. Thus it suffices to take \smash{$L \geq \frac1n
  \sum_{i=1}^n A''(x_i^\top \theta) \|x_i\|_\infty^2$}.
\begin{itemize}
\item 
For case 1: Gaussian, we have \smash{$A(u) = u^2/2$} and \smash{$A''(u) =
  1$}. Hence we may take \smash{$L = \frac1n \sum_{i=1}^n \|x_i\|_\infty^2$}. 
\item
For case 2: Bernoulli, we have \smash{$A(u) = \log(1+e^u)$} and \smash{$A''(u) =
  e^u/(1+e^u)^2 \leq \frac14$}. Hence we may take \smash{$L = \frac{1}{4n}
  \sum_{i=1}^n \|x_i\|_\infty^2$}.  
\item
For case 3: Poisson, we have \smash{$A(u)=e^u$} and \smash{$A''(u)=e^u$}. 
Again by H{\"o}lder's, we have \smash{$x_i^\top \theta \leq \|x_i\|_\infty$}
for $\theta \in \simplex$. Hence we may take \smash{$L = \frac1n \sum_{i=1}^n
  e^{\|x_i\|_\infty} \|x_i\|_\infty^2$}.  
\end{itemize}
 
\section{Proofs for Section \ref{sec:alquier}}

\subsection{Proof of Proposition \ref{prop:gibbs-equiv}}

Note \smash{$f(\theta)=\E_{\beta\sim\theta}[\hR(\beta)] = \sum_{j=1}^d
  \theta_j \hR(X_j)$}. As this is linear in $\theta$, its gradient is
constant: \smash{$\nabla_j f(\theta) = \hR(X_j)$}. The exponentiated gradient
descent iterates in \eqref{eq:egd-iterate} can therefore be written as    
\[
\theta_{t+1} (\de \beta) \propto \exp (- \eta \hR(\beta)) \cdot \theta_t (\de
\beta). 
\]
Starting at $\theta_0=u$, we see that this matches \eqref{eq:gibbs-formula} at
$1/\lambda=\eta T$. 

\subsection{Proof of Proposition \ref{prop:gibbs-risk}}

By the basic inequality \eqref{eq:md-basicineq} in Theorem
\ref{thm:md-basicineq}, 
\[
\E_{\beta\sim\thetaegdT}[\hR(\beta)] - \E_{\beta\sim\theta}[\hR(\beta)]
\leq \lambda_T\Big( \DKL(\theta, u) - \DKL(\theta, \thetaegdT) \Big).
\]
Write $\nu=\theta-\thetaegdT$, and correspondingly
\[
\E_{\beta\sim\thetaegdT}[\hR(\beta)] - \E_{\beta\sim\theta}[\hR(\beta)] = 
\int -\hR(\beta)\nu(\de\beta).
\]
Using the same representation for test risk, 
\begin{align*}
\E_{\beta\sim\thetaegdT}[R(\beta)] - \E_{\beta\sim\theta}[R(\beta)]
&= \int -R(\beta) \nu(\de\beta) \\
&= \int \Big[\hR(\beta)-R(\beta)\Big] \nu(\de\beta) 
  + \int -\hR(\beta) \nu(\de\beta) \\
&\leq \int \Big[\hR(\beta) - R(\beta)\Big]\nu(\de\beta)
  + \lambda_T \Big(\DKL(\theta, u) - \DKL(\theta,\thetaegdT)\Big) \\
&\leq \|\hR-R\|_\infty \|\nu\|_1 + \lambda_T \Big(\DKL(\theta, u) -
  \DKL(\theta,\thetaegdT)\Big) \\ 
&\leq \frac{1}{2 \lambda_T} \|\hR-R\|_\infty^2 + \frac{\lambda_T}{2} \|\nu\|_1^2
  + \lambda_T \Big(\DKL(\theta, u) - \DKL(\theta,\thetaegdT)\Big) \\ 
&\leq \frac{1}{2 \lambda_T} \|\hR-R\|_\infty^2 + \lambda_T \DKL(\theta, u).
\end{align*}
The third line uses the basic inequality, the fourth H{\"o}lder's inequality,
the fifth Young's inequality, and the last Pinsker's inequality. This proves the
first claim in the proposition.   

To prove the second claim, by a union bound and Hoeffding's inequality,
\begin{align*}
P \bigg( \|\hR-R\|_\infty > C \sqrt{\frac{\log(2d) + \delta}{2n}} \bigg)
&\leq d \cdot \sup_\beta \, \P \bigg( |\hR(\beta) - R (\beta)| > C
  \sqrt{\frac{\log(2d) + \delta}{2n}} \bigg) \\ 
&\leq 2d \exp\bigg( - \frac{2 C^2 (\log(2d) + \delta) / 2n}{C^2 / n} \bigg) 
= e^{- \delta}.
\end{align*}
Combining this with the first result in the proposition, it follows that with
probability at least \smash{$1 - e^{- \delta}$}, 
\[
\E_{\beta\sim\hthetalambda}[R(\beta)] - 
\min_{\theta: \DKL(\theta, u) \leq b} \E_{\beta \sim \theta} [R(\beta)] 
\leq \frac{C^2 (\log(2d) + \delta)}{4n \lambda} + \lambda b.
\]
Finally, setting $\lambda$ as defined in the proposition establishes the second
claim, and completes the proof. 

\section{Further details for Section \ref{sec:experiments}}

\paragraph{Optimization details: implicit regularization.}

For GD and EGD, we use learning rate schedules to cover small
$\tau$ (i.e., early training) with high resolution and to reach large $\tau$
(i.e., later training) with fewer iterations. Tables \ref{table:gd-lr-schedule}
and \ref{table:egd-lr-schedule} list the schedules for GD and EGD across 
three GLMs and two parameterization regimes. The schedule \smash{$\{(\eta^{(k)},
  T^{(k)})\}_{k=1}^m$} means that the learning rate \smash{$\eta^{(1)}$} is used
for \smash{$T^{(1)}$} iterations, then \smash{$\eta^{(2)}$} is used for the next
\smash{$T^{(2)}$} iterations, and so on.

\renewcommand{\arraystretch}{1.5}
\begin{table}[htb]
\centering\small
\begin{tabular}{| l | c | c |}
\hline
& \multicolumn{2}{c|}{\textbf{GD}} \\
\cline{2-3}
\textbf{GLM} & Underparameterized & Overparameterized \\
& $(n,d) = (200,20)$ & $(n,d) = (100, 200)$ \\
\hline
Linear & $\{(10^{-4}, 10^4), \; (10^{-3}, 10^5), \; (10^{-2}, 10^5)\}$ & same as Underparameterized \\
Logistic & same as Linear  & same as Underparameterized \\
Poisson & same as Linear & $\{(10^{-4}, 10^5), \; (2\times 10^{-4}, 2\times 10^5), \; (5\times 10^{-4}, 2\times 10^6)\}$ \\
\hline
\end{tabular}
\caption{GD learning rate schedules.}
\label{table:gd-lr-schedule}
\bigskip
    
\begin{tabular}{| l | c | c |} 
\hline
& \multicolumn{2}{c|}{\textbf{EGD}} \\
\cline{2-3}
\textbf{GLM} & Underparameterized & Overparameterized \\
& $(n,d) = (200,20)$ & $(n,d) = (30, 60)$ \\
\hline
Linear & $\{(10^{-4}, 10^5), (10^{-3}, 10^5), (10^{-2}, 10^5), (10^{-1}, 10^5)\}$ & same as Underparameterized \\
Logistic & same as Linear & same as Underparameterized \\
Poisson & same as Linear & same as Underparameterized \\
\hline
\end{tabular}
\caption{EGD learning rate schedules.}
\label{table:egd-lr-schedule}
\end{table}

\paragraph{Optimization details: explicit regularization.}

In each GLM task and in each parameterization regime, we solve $500$ regularized
problems with $\lambda$ values log-spaced on \smash{$[10^{-4}, 10^4]$}.  We use
\texttt{scipy.optimize.minimize} from the \texttt{SciPy} library; in particular,
we use the \texttt{L-BFGS-B} solver for ridge-regularized estimates and we use
the \texttt{SLSQP} solver for KL-regularized estimates.


\paragraph{Training dynamics with KL penalty.}

\begin{figure}[htb]
\centering
\includegraphics[width=0.58\textwidth]{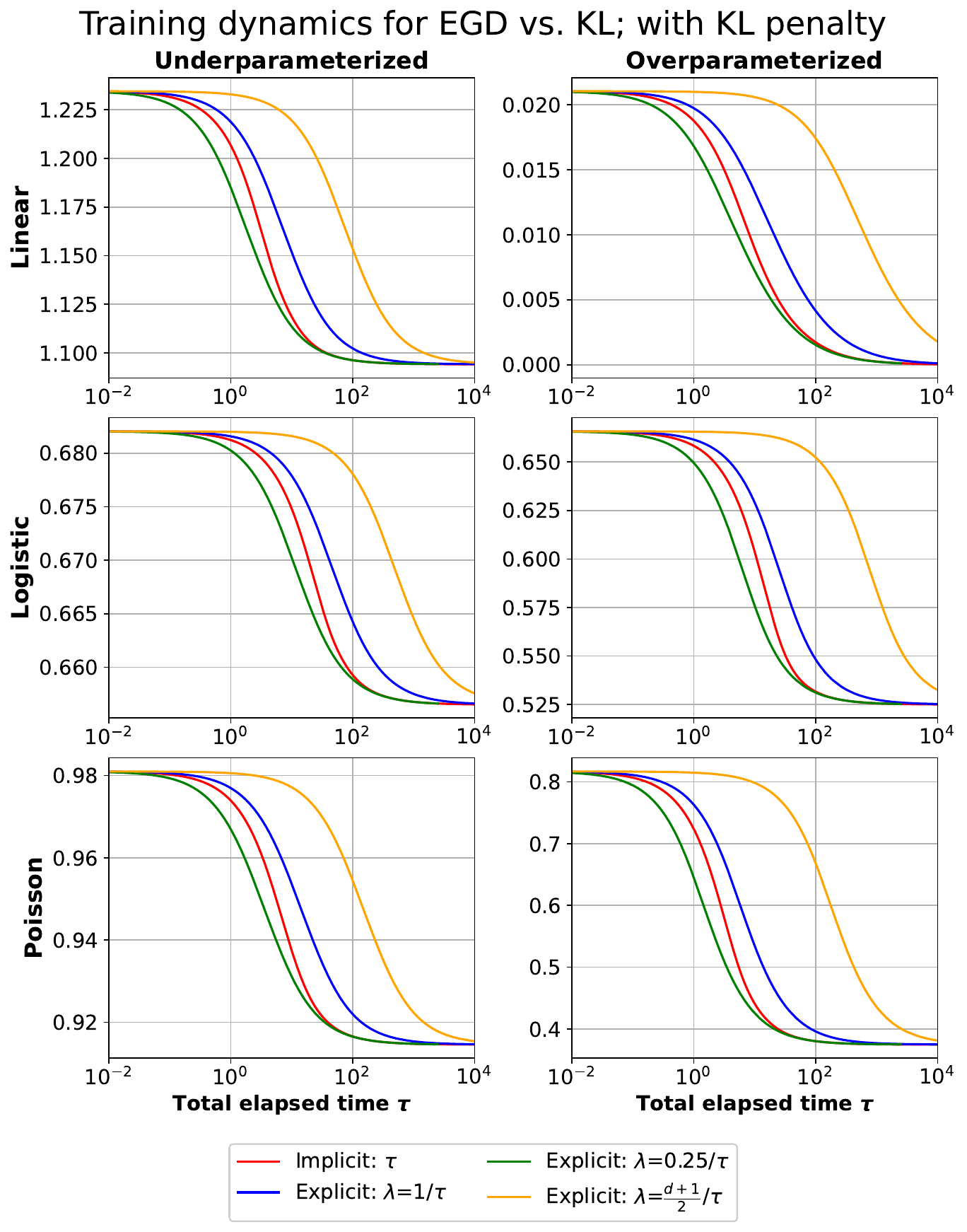}
\caption{Training envelope for model aggregation, with penalty
  $g(\theta)=\DKL(\theta,\pi)$.}  
\label{fig:training-loss-penalty-trajectory-kl}
\end{figure}

Recall in Figure \ref{fig:training-loss-penalty-trajectory} panel (b) we
investigated training dynamics for EGD vs.\ KL-regularized estimates, where
the training envelope $P(\theta,\lambda) = f(\theta) + \lambda g(\theta)$ is
defined using a penalty $g(\theta) = \|\theta-\pi\|_1^2$ (as per Corollary 
\ref{cor:md-envelope}). Figure \ref{fig:training-loss-penalty-trajectory-kl}
examines training dynamics with the envelope defined using a KL penalty. The
results are overall similar to Figure \ref{fig:training-loss-penalty-trajectory}.

\section{Proofs for Section \ref{sec:other-algos}}

\subsection{Proof of Theorem \ref{thm:proxgd-basicineq}}

The proximal gradient descent update \eqref{eq:proxgd-iterate} can be written as 
\begin{equation}
\label{eq:proxgd-iterate2}
\theta_{t+1} = \theta_t - \eta_t  G_{\eta_t}(\theta_t), 
\end{equation}
where \smash{$G_{\eta}(\theta) = \frac1\eta (\theta-\prox_{\eta
    f_1}(\theta-\eta\nabla f_0(\theta)))$} is called the generalized gradient. 
Using this notation, the proof is similar to that for the gradient descent case in 
Theorem \ref{thm:gd-basicineq}, with \smash{$G_{\eta_t}(\theta_t)$} in place of
$\nabla f(\theta_t)$. 

\medskip\noindent
\underline{\smash{Step 1:}} 
\emph{Tracking the proximity difference across adjacent iterations.} 
As before,
\[
\|\theta_t - z\|_2^2 - \|\theta_{t+1} - z\|_2^2 
= \|\theta_t - z\|_2^2 - \|\theta_t - \eta_t G_{\eta_t}(\theta_t) - z\|_2^2
= 2\eta_t \langle G_{\eta_t} (\theta_t), \theta_t - z\rangle - \eta_t^2
\|G_{\eta_t} (\theta_t)\|_2^2.  
\] 
    
\smallskip\noindent
\underline{\smash{Step 2:}}
\emph{Bounding the objective difference $f(\theta_t) - f(z)$.} 
The generalized descent lemma for proximal gradient descent (Lemma
\ref{lem:proxgd-descending} below) guarantees 
\[
2\eta_t (f(\theta_{t+1}) - f(z)) \leq 2\eta_t \langle G_{\eta_t} (\theta_t),
\theta_t - z\rangle - \eta_t^2 \|G_{\eta_t} (\theta_t)\|_2^2.
\]
Moreover, the same lemma with $z=\theta_t$ gives \smash{$f(\theta_{t+1}) \leq  
f(\theta_t) - (\eta_t / 2) \|G_{\eta_t}(\theta_t)\|_2^2$}, which in particular
shows that $f$ is nonincreasing along the iterate sequence, and
\[
2\eta_t (f(\theta_T) - f(z)) \leq 2\eta_t \langle G_{\eta_t} (\theta_t),
\theta_t - z\rangle - \eta_t^2 \|G_{\eta_t} (\theta_t)\|_2^2.
\]
Combined with the result of Step 1, 
\[
2\eta_t(f(\theta_T) - f(z)) \leq \|\theta_t - z\|_2^2 - \|\theta_{t+1} -
z\|_2^2. 
\]
The remainder of the proof (Step 3) follows exactly as before. 

\begin{lemma}
\label{lem:proxgd-descending}
Let $f_0:\R^d\to \R$ be convex, differentiable, and $L$-smooth with $L>0$, and 
let $f_1:\R^d\to \R$ be convex. For $f=f_0+f_1$ and any $z \in \R^d$, the
proximal gradient update in \eqref{eq:proxgd-iterate2} with \smash{$\eta_t
  \in(0,2/L]$} satisfies  
\[
f(\theta_{t+1}) \leq f(z) + \langle G_{\eta_t}(\theta_t), \theta_t - z
\rangle - \eta_t\Big(1-\frac L2 \eta_t\Big) \|G_{\eta_t}(\theta_t)\|_2^2.
\]
For $\eta_t \leq 1/L$, this implies
\[
f(\theta_{t+1}) \leq f(z) + \langle G_{\eta_t}(\theta_t), \theta_t - z
\rangle - \frac{\eta_t}{2} \|G_{\eta_t}(\theta_t)\|_2^2.
\]
\end{lemma}

\begin{proof}
By $L$-smoothness,
\begin{align*}
f_0(\theta_{t+1}) &\leq f_0(\theta_t) + \langle \nabla f_0(\theta_t),
\theta_{t+1}-\theta_t \rangle + \frac{L}{2} \|\theta_{t+1}-\theta_t\|_2^2 \\ 
&\leq f_0(\theta_t) - \langle \nabla f_0(\theta_t), \eta_t G_{\eta_t}(\theta_t)
\rangle + \frac{L \eta_t^2}{2} \|G_{\eta_t}(\theta_t)\|_2^2 \\
&\leq f_0(z) + \langle \nabla f_0(\theta_t), \theta_t - 
\eta_t G_{\eta_t}(\theta_t) - z \rangle + \frac{L \eta_t^2}{2}
\|G_{\eta_t}(\theta_t)\|_2^2, 
\end{align*}
where in the last line we used $f_0(\theta_t) \leq f_0(z) + \nabla
f_0(\theta_t)^\T (\theta_t - z)$ by convexity. Moreover, by the subgradient 
optimality condition defining the proximal operator,    
\[
\theta_t - \eta_t \nabla f_0(\theta_t) - \theta_{t+1} \in \eta_t \partial
f_1(\theta_{t+1}),
\]
which may be equivalently written as
\[
G_{\eta_t}(\theta_t) - \nabla f_0(\theta_t) \in \partial f_1(\theta_{t+1}).
\]
This means that 
\[
f_1(\theta_{t+1}) \leq f_1(z) + \langle G_{\eta_t}(\theta_t) - \nabla
f_0(\theta_t), \theta_t - \eta_t G_{\eta_t}(\theta_t) - z \rangle. 
\]
Combining this with the above inequality on $f_0$, and using $f=f_0+f_1$, we
obtain
\[
f(\theta_{t+1}) \leq f(z) + \langle G_{\eta_t}(\theta_t), \theta_t - z
\rangle - \eta_t \|G_{\eta_t}(\theta_t)\|_2^2 + \frac{L \eta_t^2}{2}
\|G_{\eta_t}(\theta_t)\|_2^2, 
\]
which can be regrouped to give the desired result.
\end{proof}

\subsection{Proof of Theorem \ref{thm:nolips-basicineq}}

The proof is similar to that for the mirror descent case in Theorem
\ref{thm:md-basicineq}. 

\medskip\noindent
\underline{\smash{Step 1:}} 
\emph{Tracking the proximity difference across adjacent iterations.} 
By the same argument as before,  
\[
\eta_t \langle \nabla f(\theta_t), \theta_{t+1} - z\rangle
\leq \Dphi(z,\theta_t) - \Dphi(z, \theta_{t+1}) - \Dphi(\theta_{t+1}, \theta_t).
\]

\smallskip\noindent
\underline{\smash{Step 2:}}
\emph{Bounding the objective difference $f(\theta_t) - f(z)$.} 
By convexity of $L\phi-f$, 
\[
L\phi(\theta_{t+1}) - f(\theta_{t+1}) \geq L\phi(\theta_t) - f(\theta_t) + 
\langle L\nabla \phi(\theta_t) - \nabla f(\theta_t), \theta_{t+1}-\theta_t
\rangle, 
\]
and rearranging gives 
\[
f(\theta_{t+1}) \leq f(\theta_t) + \langle \nabla f(\theta_t),
\theta_{t+1}-\theta_t \rangle + L D_\phi(\theta_{t+1},\theta_t). 
\]
Moreover, $f(\theta_t) \leq f(z) + \langle \nabla f(\theta_t), \theta_t - z
\rangle$ by convexity, thus 
\[
f(\theta_{t+1}) \leq f(z) + \langle \nabla f(\theta_t),
\theta_{t+1}-z \rangle + L D_\phi(\theta_{t+1},\theta_t). 
\]
Combining this with the result of Step 1, we have 
\begin{align*}
\eta_t (f(\theta_{t+1}) - f(z))
&\leq \Dphi(z,\theta_t) - \Dphi(z, \theta_{t+1}) - (1-L\eta_t)
\Dphi(\theta_{t+1}, \theta_t) \\
&\leq \Dphi(z,\theta_t) - \Dphi(z, \theta_{t+1}),
\end{align*}
where we used $\eta_t\leq 1/L$ in the last inequality. Taking $z=\theta_t$, this
shows us that $f$ is nonincreasing along the iterate sequence, and therefore
\[
\eta_t (f(\theta_T) - f(z)) \leq \Dphi(z,\theta_t) - \Dphi(z, \theta_{t+1}).
\]
The remainder of the proof (Step 3) follows exactly as before. 

\end{document}